\documentclass[a4paper,12pt,twoside]{article}
\usepackage{graphicx}
\usepackage{amsmath,amssymb,amsthm,textcomp}
\usepackage[english]{babel}

\newcounter{ass}

\RequirePackage{amsmath} \RequirePackage{amssymb}
\def\eps{\varepsilon}

\def\N{\mathbf{N}}
\def\Z{\mathbf{Z}}
\def\R{\mathbb{R}}

\def\P{\mathbf{P}}
\def\E{\mathbf{E}}
\def\F{\mathcal{F}}

\newcommand{\bpf}  {\noindent{\sc Proof} : }
\newcommand{\epf}  {\hfill $\square$}
\newtheorem{theorem}{Theorem}[section]

\newtheorem{proposition}[theorem]{Proposition}
\newtheorem{lemma}[theorem]{Lemma}
\newtheorem{corollary}[theorem]{Corollary}
\newtheorem{remark}[theorem]{Remark}

\numberwithin{equation}{section}
\textheight=8.5in \numberwithin{equation}{section}
\begin{document}
\title{The $\Lambda$-lookdown model with selection}
\author{B. Bah and E. Pardoux}
\date{\today}
\maketitle

\abstract{ The goal of this paper is to study the lookdown model with selection  in the case of a population containing two types of individuals, with a reproduction model which is dual to the $\Lambda$-coalescent. In particular we formulate the infinite population ``$\Lambda$-lookdown model with selection". When the measure $\Lambda$ gives no mass to $0$, we show that the proportion of one of the two types converges, as the population size $N$ tends to infinity, towards the solution to a stochastic differential equation driven by a Poisson point process. We show that one of the two types fixates in finite time if and only if the $\Lambda$-coalescent comes down from infinity. 
We give precise asymptotic results in the case of the Bolthausen--Sznitman coalescent. We also consider the general case of a combination of the Kingman and the $\Lambda$-lookdown model.}


\section{Introduction}
In this paper we consider the lookdown (which is in fact usually called the ``modified lookdown'') model with selection where we replace the usual reproduction model by a population model dual to the $\Lambda$-coalescent. We first recall the models from \cite{Pitman} and \cite {DK96}, and then we will describe the variant which will be the subject of the
present paper.

Pitman \cite{Pitman} and Sagitov \cite{Sagitov} have pointed at an important class of exchangeable coalescents  whose laws can be characterized by an arbitrary finite measure $\Lambda$ on [0, 1]. Specifically, a $\Lambda$-coalescent is a Markov process ($\Pi_t, t\ge 0$) on ${\mathcal P}_{\infty}$ (the set of partition of $\N$) started from the partition $0_{\infty}:=\{\{1\}, \{2\}, \dots\}$ and such that, for each integer $n\ge 2$, its restriction $(\Pi^{[n]}_t, t\ge 0)$ to ${\mathcal P}_n$ (the set of partitions of $\{1, 2, \dots, n\}$) is a continuous time Markov chain that evolves by coalescence events, and whose evolution can be  described as follows.\\
Consider the rates
\begin{equation}\label{tau}
\lambda_{k, \ell}=\int_{0}^{1}p^{\ell-2}(1-p)^{k-\ell}\Lambda(dp),\quad \mbox{$2\le \ell \le k$.}
\end{equation}
Starting from a partition in ${\mathcal P}_n$ with $k$ non-empty blocks, for each $\ell=2,\dots,k,$ every possible merging of $\ell$ blocks (the other $k-\ell$ blocks remaining unchanged) occurs at rate $\lambda_{k, \ell}$, and no other transition is possible. This description of the restricted processes $\Pi^n$ determines the law of the $\Lambda$-coalescent $\Pi$.\\
Note that if $\Lambda(\{0\})=\Lambda([0, 1])>0$, then only pairwise merging occurs, and the corresponding $\Lambda$-coalescent is just a time rescaling (by $\Lambda({0})$) of the Kingman coalescent. When $\Lambda(\{0\})=0$ which we will assume except in the very last section of this paper, a realization of the $\Lambda$-coalescent can be constructed (as in \cite{Pitman}) using a Poisson point process 
\begin{equation}\label{measure}
m=\sum_{i=1}^{\infty}\delta_{t_i, p_i}
\end{equation}
on $\R_{+} \times (0, 1]$ with intensity measure $dt \otimes \nu(dp)$ where $\nu(dp)=p^{-2}\Lambda(dp)$. We will assume that the measure $\nu(dp)$  has infinite total mass. Each atom $(t, p)$ of $m$  influences the evolution as follows :
\begin{itemize}
\item for each block of $\Pi(t^-)$ run an independent Bernoulli ($p$) random variable; 
\item all the blocks for which the Bernoulli outcome equals  1 merge immediately 

into one single block, while all the other blocks remain unchanged.   
\end{itemize}
In order to obtain a construction for a general measure $\Lambda$, one can superimpose onto the $\Lambda$-coalescent independent pairwise mergers at rate $\Lambda(\{0\})$.

The lookdown construction was first introduced by Donnelly and Kurtz in 1996 \cite{DK96}. Their goal was to give a construction of the Fleming-Viot superprocess that provides an explicit description of the genealogy of the individuals in a population. Donnelly and Kurtz subsequently modified their construction in \cite{DK99a} to include more general measure-valued processes. Those authors extended their construction to the selective and recombination case \cite{DK99b}.   

We are going to present our model which we call $\Lambda$-lookdown model with selection. An important feature of our model is that we will describe it for a population of infinite size, thus retaining the great power of the lookdown construction. As far as we know, this has not yet been done in the case of models {\it with selection} except in our previous publication \cite{bb.ep.bs}, where we considered a model dual to Kingman's coalescent.  

We consider the case of two alleles $b$ and $B$, where $B$ has a selective advantage over $b$. This selective advantage is modelled by a death rate $\alpha$ for the type $b$ individuals.  We will consider the proportion of $b$ individuals. The type $b$ individuals are coded by 1, and the type $B$ individuals by 0. We assume that the individuals are placed at time $0$ on levels $1, 2, \dots,$ each one being, independently from the others, 1 with probability $x$, 0 with probability $1-x$, for some $0<x<1$. For each $i\ge 1$ and $t\ge 0$, let $\eta_t(i) \in \{0, 1\}$ denote the type of the individual sitting on level $i$ at time $t$. The evolution of $(\eta_t(i))_{i\ge1}$ is governed by the two following mechanisms.
\begin{enumerate}
\item{\it{Births}}
Each atom $(t, p)$ of the Poisson point process $m$ corresponds to a birth event. To each $(t, p)\in m$, we associate a sequence of i.i.d Bernoulli random variables $(Z_i, i\ge 1)$ with parameter $p$. Let
$$I_{t, p}=\{i\ge 1 : Z_i=1\}.$$
and
$$\ell_{t, p}=\inf\{i\in I_{t, p} : i>\mbox{min\;} I_{t, p}\}$$
At time $t$, those levels with $Z_i$=1 and $i\ge \ell_{t, p}$ modify their label to $\eta_{t^-}(\mbox{min\;} I_{t, p})$. In other words, each level in $I_{t, p}$ immediately adopts the type of the smallest level participating in this birth event. For the remaining levels, we reassign the types so that their relative order immediately prior to this birth event is preserved. More precisely    

\begin{eqnarray*}
\eta_t(i)=
\begin{cases}
\eta_{t^-}(i), &  \mbox{if} \;  \;    i<\ell_{t, p}  \\
\eta_{t^-}(\mbox{min\;}I_{t, p}), &  \mbox{if} \;  \;   i\in I_{t, p}\setminus\{\mbox{min}\;I_{t, p}\} \\
\eta_{t^-}(i-(\#\{I_{t, p}\cap [1,\dots,i]\}-1)), &  \mbox{}      \mbox{otherwise}
\end{cases}
\end{eqnarray*}
We refer to the set $I_{t, p}$ as a multi-arrow at time $t$, originating from min $I_{t, p}$, and with tips at all other points of $I_{t, p}.$
This procedure is usually referred to as the modified lookdown construction of Donnelly and Kurtz. In the original construction, the types of the levels in the complement of $I_{t, p}$ remained unchanged at time $t$, hence the types $\eta_{t^-}(i)$, for $i\in I_{t, p}\setminus\{\mbox{min}\;I_{t, p}\}$ got erased from the population at time $t$.

\item{\it Deaths} Any type 1 individual dies at rate $\alpha$, his vacant level being occupied by his right neighbor, who himself is replaced by his right neighbor, etc. In other words, independently of the above arrows, crosses are placed on all levels according to mutually independent rate $\alpha$ Poisson processes. Suppose there is a cross at level $i$ at time $t$. If $\eta_{t^-}(i)=0$,
nothing happens. If $\eta_{t^-}(i)=1$, then \begin{eqnarray*}
\eta_t(k)=
\begin{cases}
\eta_{t^-}(k), &  \mbox{if} \;  \;    k<i;  \\
\eta_{t^-}(k+1), &  \mbox{if} \;  \;   k\ge i .\\
\end{cases}
\end{eqnarray*}
\end{enumerate}
We refer the reader to Figure 1 for a pictural representation of our model. Note that the type of the newborn individuals are found by ``looking down'', 
while the type of the individual who replaces a dead individual is found by looking up. So maybe our model could be called ``look-down, look-up''.  

\begin{figure}
\begin{center}
\scalebox{0.265}{
\includegraphics{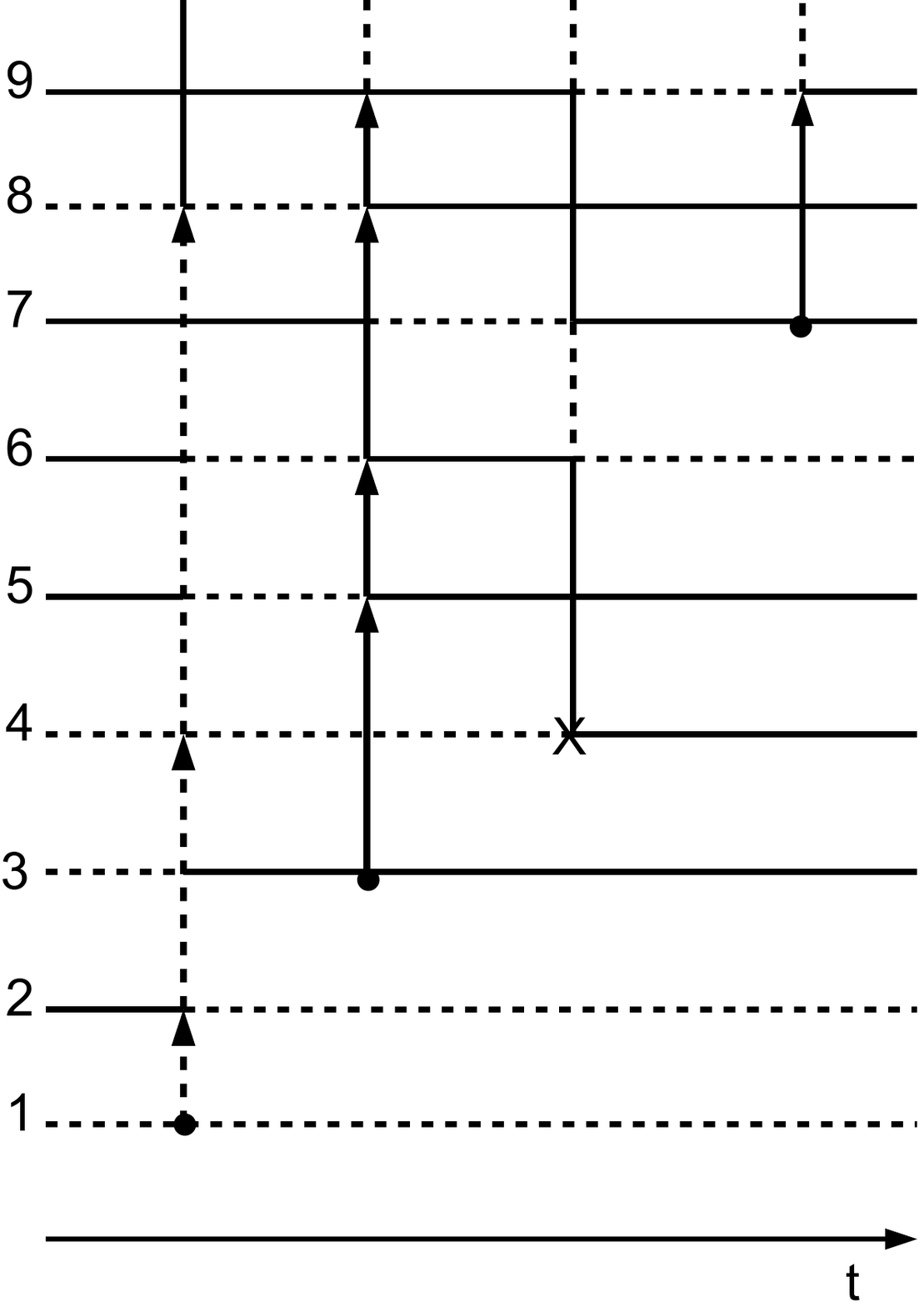}
} 
\caption{The graphical representation of the $\Lambda$-lookdown model with selection of size $N=9$. Solid lines represent type $B$ individuals, while dotted lines represent type $b$ individuals. 
}
\label{figure1}
\end{center}
\end{figure}
Since we have modelled selection by death events, the evolution of the $N$ first individuals $\eta_t(1),\dots,\eta_t(N)$ depends upon the next ones, and $X^N_t=N^{-1}(\eta_t(1)+\dots+\eta_t(N))$, the proportion of type $b$ individuals among the $N$ lowest levels,  is not a Markov process. We will show however that for each $t>0$ the collection of r.v.'s $\{\eta_t(k),\ k\ge1\}$ is well defined (which is not obvious in our setup) and constitutes an exchangeable sequence of $\{0,1\}$--valued random variables. We can then apply de Finetti's theorem, and prove that $X^N_t\to X_t$ a.s for any fixed $t\ge 0$, where $(X_t)_{t\ge0}$ is a $[0,1]$--valued Markov process, which is a solution to the stochastic differential equation (which we call the $\Lambda$--Wright--Fisher SDE
with selection)
\begin{equation}\label{eqXt}
X_t = x- \alpha \int_{0}^{t} X_s(1-X_s)ds+\int_{[0, t]\times]0, 1[^2}p({\bf 1}_{u\le X_{s^-}}-X_{s^-})\bar{M}(ds, du, dp),
\end{equation}
where $\bar{M}(ds, du, dp)=M(ds, du, dp)-p^{-2}dsdu\Lambda(dp)$, and $M$ is a Poisson point measure on $\R_{+} \times ]0, 1[\times]0, 1]$ with intensity $dsdup^{-2}\Lambda(dp)$. 
The process $(X_t)_{t\ge0}$ represents the proportion of type $b$ individuals at time $t$ in the infinite size population. Note that uniqueness of a solution to \eqref{eqXt} is proved in \cite{DL}.  

The paper is organized as follows. We both construct our process, and establish the crucial exchangeability property satisfied by the $\Lambda$-lookdown model with selection in section 2.  In section 3 we establish the convergence of $X^N$ to the solution to \eqref{eqXt}.  In section 4 we show that one of the two types fixates in finite time if and only if the $\Lambda$-coalescent comes down from infinity. Moreover, in the case of no fixation, we show that $X_t \rightarrow X_\infty \in \{0, 1\}$ as $t\rightarrow \infty$, and 
discuss when $X_\infty=0\ a.s$ and when $\P(X_{\infty}=1)>0$. 
In the case of the Bolthausen--Sznitman coalescent (which does not come down from infinity), we precise the law of $X_\infty$, and study the speed at which either 
of the two types invades the whole population.
Finally, we extend our results to the case $\Lambda(\{0\})>0$ in the last section 5.   

In this paper, we use $\N$ to denote the set of positive integers $\{1,2,\ldots\}$, and $[n]$ to denote the set $\{1,\dots,n\}$. 
We suppose that the measure $\Lambda$ fulfills the condition
\begin{equation}\label{condition}
0<\Lambda((0, 1))<\infty, \quad \Lambda(\{1\})=0, 
\end{equation}
and in all the paper except in section 5, we assume that $\Lambda(\{0\})=0$.

\section{The lookdown process, exchangeability}\label{sectLD}
\subsection{Some results for general $\Lambda$ }
Throughout the paper, the notation
$$\mu_r:=\int_{[0, 1]}p^r \Lambda(dp)$$
is used for the $r$th moment of the finite measure $\Lambda$ on [0, 1] for arbitrary real $r$. Note that $\mu_r$ is a decreasing function of $r$ with $\infty >\mu_0\ge\mu_r > 0$ for $r\ge 0$, while $\mu_r$ may be either finite or infinite for $r<0$.
For $r = 0, 1, \cdots$observe from
\eqref{tau} that $\mu_r = \lambda_{r+2, r+2}$ is the rate at which $\Pi_n$ jumps to its absorbing state
\{[n]\} from any state with $r+2$ blocks. Let $X$ denote a random variable with distribution $\mu_0^{-1}\Lambda$, defined on
some background probability space $(\Omega, \mathcal {F}, \P) $ with expectation operator $\E$, so
$\E(X^r) = \mu_r/\mu_0$. Recall the formula \eqref{tau} for the transition rates $\lambda_{k, \ell}$ of the $\Lambda$-coalescent,
which we rewrite as
$$\lambda_{k, \ell}=\mu_0\E(X^{\ell-2}(1-X)^{k-\ell})\;\; \mbox{for all}\; 2\le \ell\le k. $$

For any partition with a finite number $n\ge 2$ of blocks, the total rate of transitions of all kinds in a $\Lambda$-coalescent, which can be rewritten as
\begin{align*}
\lambda_n:&=\sum_{\ell=2}^{n}{n\choose \ell}\lambda_{n, \ell}=\int_{0}^{1}\frac{1-(1-p)^n-np(1-p)^{n-1}}{p^2}\Lambda(dp)\label{total-rate}\\
&=\mu_0\E\left[\frac{1-(1-X)^n-nX(1-X)^{n-1}}{X^2}\right].
\end{align*}
By monotone convergence, 

\begin{equation*}
\lambda_n \uparrow \mu_{-2}=\int_{[0, 1]}p^{-2}\Lambda(dp)\;\; \mbox{as} \;\; n \uparrow \infty.
\end{equation*}
\subsection{Construction of our process}
In this section, we will construct the process $\{\eta_t(i), i\ge 1, t\ge 0\}$ corresponding to a given initial condition $(\eta_0(i), i\ge 1)$ defined in the Introduction. 

Recall the Poisson point process $m$ defined in \eqref{measure}. For each $n\ge 1$ and $ t\ge 0$ , let
$$I(n, t)=\{k\ge 1: t_k\in [0, t]\; \text{and} \;\#\{I_{t_k, p_k}\cap [n]\}\ge 2\}.$$
We have
\begin{lemma}\label{evolution}
For each $n\ge 1$ and $t\ge 0$,
$$\# I(n, t)<\infty\quad \text{a.s.} $$
\end{lemma}
\bpf
Each atom $(t, p)$ of $m$ affects at least 2 of the $n$ first individuals with probability 
$$1-(1-p)^n-np(1-p)^{n-1}\le {n \choose 2}p^{2}.$$
Consequently 
$$\E(\# I(n, t))\le {n \choose 2}t\int_{0}^{1}\Lambda(dp)<\infty.$$
The result follows.
\epf

\subsubsection{$\Lambda$-lookdown model without selection}

In this subsection, we essentially follow \cite{DK99a}.
For each $N\ge 1$, one can define the vector $\xi^N_t=(\xi^N_t(1),\dots,\xi^N_t(N)), t\ge 0$ with values in $\{0, 1\}^N$, by 
\begin{enumerate}
\item $\xi^N_0(i):=\eta_0(i)\; \mbox{for all}\ i\ge 1$.
\item At any birth event $(t, p) \in m$ and such that $\{I_{t, p}\cap [N]\}\ge 2$, for each $i\in [N]$, $\xi^N_t(i)$ evolves as follows 
\begin{eqnarray*}
\xi^N_t(i)=
\begin{cases}
\xi^N_{t^-}(i), &  \mbox{if} \;  \;    i<\ell_{t, p}  \\
\xi^N_{t^-}(\mbox{min\;}I_{t, p}), &  \mbox{if} \;  \;   i\in I_{t, p}\setminus\{\mbox{min}\;I_{t, p}\} \\
\xi^N_{t^-}(i-(\#\{I_{t, p}\cap [1,\dots,i]\}-1)), &  \mbox{}      \mbox{otherwise}.
\end{cases}
\end{eqnarray*}
\end{enumerate}
Using the above lemma, we see that the process $\xi^N_t$ has finitely many jumps on $[0, t]$ for all $t>0$, hence its evolution is well defined. From this definition, one can easily deduce that the evolution of the type at levels 1 up to $i$ depends only upon the types at levels up to $i$. Consequently, if $1\le N<M$, the restriction of $\xi^M$ to the $N$ first levels yields $\xi^N$, in other words :
$$\{\xi^M_t(1), \dots, \xi^M_t(N), t\ge 0\} \equiv \{\xi^N_t(1), \dots, \xi^N_t(N), t\ge 0\}.$$  
Hence, the process $\eta=\xi^{\infty}$ is easily defined by a projective limit argument as a $\{0, 1\}^{\infty}$-valued process.
\subsubsection{$\Lambda$-lookdown model with selection}\label{Llms}
This section is devoted to the construction of the infinite population lookdown model with selection. 

For each $M\ge 1$, we consider the process $(\eta^M_t(i), i\ge 1, t\ge 0)$ obtained by applying all the arrows between $1\le i<j<\infty$, and only the crosses on levels 1 to $M$. Using the fact that we have a finite number of crosses on any finite time interval, it is not hard to see that the process $(\eta^M_t, t\ge 0)$ is well defined by applying the model without selection between two consecutive crosses, and applying the recipe described in the Introduction at a death time. More generally, our model is well defined if we suppress all the crosses above a curve which is bounded on any time interval $[0, T]$. Note also that, if we remove or modify the arrows and or the crosses above the evolution curve of a type $B$ individual, this does not affect her evolution as well as that of those sitting below her.  

At any time $t\ge 0$, let $K_t$ denote the lowest level occupied by a $B$ individual. Of course, if $K_0=1$, then $K_t=1$, for all $t\ge 0$. If for any $T$, $\sup_{0\le t\le T} K_t<\infty$ a.s, then the process $\{K_t, t\ge 0\}$ is well defined by taking into account only those crosses below the curve $K_t$, and evolves as follows. When in state $n>1$, $K_t$ jumps to 
\begin{enumerate}
\item $n+k$ at rate  ${n+k-1\choose k+1} \lambda_{n+k, k+1} , \quad k\ge 1 $; 
\item $n-1$  at rate $\alpha (n-1)$, $\alpha>0$,
\end{enumerate} 
where we have used the notation $\lambda_{k,\ell}$ defined by \eqref{tau}. In other words,
the infinitesimal generator of the Markov process $\{K_t, t\geq 0\}$ is given by:
\begin{equation}\label{generator}
\mathcal{L}g(n)=\sum_{k=1}^{\infty} {n+k-1\choose k+1} \lambda_{n+k, k+1}[g(n+k)-g(n)]+\alpha (n-1)[g(n-1)-g(n)].
\end{equation}

Now, we are going to show that the process $\{\eta_t(i), i\ge 1, t\ge 0\}$ is well defined. For this, we study two cases.

\textbf{Case 1:} $K_t \rightarrow \infty$ as $t\rightarrow \infty.$

For each $N\ge 1, t\ge 0$, we define
$$K^N_t= \mbox{the level of the $N-$ th individual of type $B$ at time $t$};$$
and 
$$T^{N}_{\infty}=\inf\{t\ge 0 : K^N_t=\infty\}.$$
We have  $T^1_\infty\ge T^2_\infty\ge \cdots> 0$. For each $N\ge 1$, we define
$$\mathcal{H}_N=\{(s, k); k\le K^N_s\}.$$
Consider first the event
$$A= \{ T^N_\infty=\infty\; ,\forall \  N\ge 1\}\footnote{We shall see below that $\P(A)=1$ if the $\Lambda$-coalescent does not come down from infinity and $\P(A)=
\P(\text{type B fixates})$ otherwise.}.$$
Recall the Poisson point measure $m$ defined in \eqref{measure}. Now, for each $N\ge 1$, we define the process$(\eta^N_t(i), i\ge 1, t\ge 0)$, with values in $\{0, 1\}^{\infty}$, by
\begin{enumerate}
\item $\eta^N_0(i):=\eta_0(i)\; \mbox{for all}\ i\ge 1$.
\item At any birth event $(t, p)\in m, \eta^N_t$ evolves as follows 
\begin{eqnarray*}
\eta^N_t(i)=
\begin{cases}
\eta^N_{t^-}(i), &  \mbox{if} \;  \;    i<\ell_{t, p}  \\
\eta^N_{t^-}(\mbox{min\;}I_{t, p}), &  \mbox{if} \;  \;   i\in I_{t, p}\setminus\{\mbox{min}\;I_{t, p}\} \\
\eta^N_{t^-}(i-(\#\{I_{t, p}\cap [1,\dots,i]\}-1)), &  \mbox{}      \mbox{otherwise},
\end{cases}
\end{eqnarray*}

\item Suppose there is a cross on level $j$ at time $s$. If $(s, j) \notin \mathcal{H}_N$ or $(s, j) \in \mathcal{H}_N$ and $\eta_{s^-}(j)=0$, nothing happens. If $(s, j)\in \mathcal{H}_N$ and $\eta_{s^-}(j)=1$, then
\begin{eqnarray*}
\eta^N_s(i)=
\begin{cases}
\eta^N_{s^-}(i), &  \mbox{if} \;  \;    i<j;  \\
\eta^N_{s^-}(i+1), &  \mbox{if} \;  \;   i\ge j.\\
\end{cases}
\end{eqnarray*}
\end{enumerate} 
In other words, the process $\{\eta^N_t(i), i\ge 1, t\ge 0\}$ is obtained by applying all the arrows between $1\le i<j<\infty$, and only the crosses on levels 1 to $K^N_t$. On the event $A$, we have a finite number of such crosses on any finite time interval, and $(\eta^N_t(i), i\ge 1, t\ge 0)$ is constructed  as explained above. Now, let
$$\mathcal{H}=\cup_N \mathcal{H}_N.$$
By a projective limit argument, we can easily deduce that the process $\{\eta_t(i), i\ge 1, t\ge 0\}$ is well defined on the set $\mathcal{H}$. Our model is defined on the event $A$.

Now we consider the event $A^{c}$. We first work on the event $\{T^1_{\infty}<\infty\}$. This means that the allele $b$ fixates in finite time. It implies that for each $N\ge 2, T^N_\infty$ is finite as well. Consider first the process $\{\eta^1_t(i), i\ge 1, t\ge 0 \}$ 
defined on $\mathcal{H}_1$, i.e we take into account all the arrows between $1\le i<j\le K^1_t$, and only the crosses on levels $1$ to $K^1_t$. This process is well defined on the time interval $[0, T^1_{\infty})$. However, on the interval $[T^1_{\infty}, \infty), \eta^1_t(i)=1, \forall i\ge 1,$ hence the process is well defined in $\mathcal{H}_1$. We next consider the process $\{\eta^2_t(i), i\ge 1, t\ge 0\}$ 
defined on $\mathcal{H}_2$. This process is well defined on the time interval $[0, T^2_{\infty})$. But on the interval $[T^2_{\infty}, \infty)$, there is at most one $B$, whose position is completely specified from the previous step. Iterating that procedure, and using again a projective limit argument, we define the full $\Lambda$-lookdown model with selection.

If $T^1_{\infty}=+\infty$, but $T^N_{\infty}<+\infty$ for some $N$, the construction is easily adapted to that case. In fact some arguments in section 4 below show that this cannot happen with positive probability.\\\\
\textbf{Case 2 :} $K_t \nrightarrow \infty, t\rightarrow \infty .$

Let 
$$T_1=\inf\{t\ge 0 : K_t=1\}.$$
We now show that $\{T_1<\infty\}$ a.s. on the set $\{K_t \nrightarrow \infty, t\rightarrow \infty\}$.
Indeed, for any stopping time $T$ and $M>1$, define 
$D_{T,M}$ to be the event that there is at least one cross on each of the levels $1,2,\ldots,M-1$ on the the interval $(T,T+1)$, and 
$B_{T,M}$ to be the event that no birth arrow points to a level less than or equal to $M$ on the time interval $(T,T+1)$. It is plain that the quantity
$$p_{\alpha,M}=\P(D_{T,M}\cap B_{T,M}|\F_T)$$ is deterministic, independent of $T$, and that $p_{\alpha,M}>0$. Now clearly
$$\{K_T\le M\}\cap D_{T,M}\cap B_{T,M}\subset \{K_{T+1}=1\}.$$
Hence
$$\P(K_{T+1}=1|K_T\le M)\ge p_{\alpha,M},$$
or equivalently
$$\P(K_{T+1}>1|K_T\le M)\le 1-p_{\alpha,M}.$$
Let now 
\begin{equation*}
A_M:=\left\{\begin{aligned}
&\text{ there exists an infinite sequence of stopping times }T_M^k\\&\text{such that }
T_M^{k+1}\ge T_M^k+1\text{ and }K_{T_M^k}\le M,\text{ for all } k\ge1.
\end{aligned}
\right\}
\end{equation*}
 We deduce from the last inequality and  the strong Markov property that for any $n\ge1$,
$$\P(A_M\cap\{K_{T_M^{n+1}}>1\})\le (1-p_{\alpha,M})^n.$$
consequently $\P(A_M\cap\{T_1=+\infty\})=0$. This being true for all $M>0$, the claim follows.

If $T_1<\infty$,  the idea is to show that there exists an increasing mapping $\psi : \N\rightarrow \N$ such that a.s. for $N$ large enough, any individual sitting on level $\psi(N)$ at any time never visits a level below $N$, with the convention that if that individual dies, we replace him by his neighbor below. Once this is true, the evolution of the individuals sitting on levels $1, 2, \dots, N$ is not affected by deleting the crosses above level $\psi(N)$. Hence it is well defined. If this holds for all $N$ large enough, the whole model is well defined.

Let
$$M=\sup_{0\le t< T_1}K_t.$$
For each $N\ge M$, we will show that an individual sitting on a high enough level at any time $t\ge 0$ never visits a level below $N$. In order to prove this, we couple our model with the following one. 

On the interval $[0, T_1]$, we erase all the arrows pointing to levels above $K_t$, and pretend that all individuals above level $K_t$, $0\le t\le T_1$, are of type $b$, i.e coded by 1, and we apply all the crosses above level $K_t$. This model is clearly well defined since until $T_1$ there is only one $0$, all other sites being occupied by 1's. We next extend this model for $t>T_1$ as follows :

For each $t\ge T_1$, let $\bar{K_t}$ denote the lowest level occupied by a $b$ individual. At time 
$T_1$, $\eta_{T_1}(1)=0, \eta_{T_1}(i)=1,$ for all $i\ge 2$. At any time $t>T_1$, we shall have $\eta_t(i)=0$ 
for $i< \bar{K_t}$, and $\eta_t(i)=1$ for $i\ge\bar{K_t}$. Again all crosses are kept, and we keep only 
those arrows whose tip hits a level $j\le \bar{K_t}$.

This model is well defined. For each $N\ge 1$, we define $S_N$ as the first time where all the $N$ first individuals of this model are of type $B$. We have
\begin{lemma}\label{Nfixate}
If $T_1<\infty$, then for each $N\ge 1$,
$$S_N<\infty \ a.s$$
\end{lemma}
\bpf
The result follows from $T_1<\infty$ and the fact that the process of arrows from 1 to 2 is a Poisson process with rate $\lambda_{2,2}=\Lambda((0, 1))$.

\epf
\bigskip

Now, let $\varphi(N)=Ne^{\alpha S_N}(Ne^{\alpha S_N}+1)+K_0$ and $\{\xi^{\varphi(N)}_t,\ t\ge0\}$ denote the process which describes the position at time $t$ of the individual sitting on level $\varphi(N)$ at time $0$ in the present model.

We will prove below that
the individual who sits on level $\varphi(N)$ at time 0 will remain below the level $\varphi(N)+N$ on the time interval $[0, S_N]$. If she does not visit any level below $N$ before time  
$S_N$, she will never visit any level below $N$ at any time, and moreover any individual who visits level $\varphi(N)+N$ before time $S_N$ will remain above the individual who was sitting at level $\varphi(N)$ at time 0 until $S_N$, hence will never visit any level below $N$.

Since the ``true" model has more arrows and less ``active crosses" than the present model, if we show that in the present model a.s. there exists $N$ such that the individual who starts from level $\varphi(N)$ at time 0 never visits a level below $N$, we will have that in the true model a.s. for $N$ large enough the evolution within the box $(t, i)\in [0, \infty)\times \{1, 2, \dots, N\}$ is not altered by removing all the crosses above $\varphi(N)+N$. A projective limiting argument allows us then to conclude that the full model is well defined.   

The result will follow from the Borel-Cantelli lemma and the following lemma.
 \begin{lemma}\label{MajorD}
 If $T_1< \infty$, then for each $N\ge M$,
$$\widehat{\P}_N(\exists 0<t\le S_N\ \text{such that } \xi^{\varphi(N)}_t\le N)\le \frac{2}{N^2},
$$
where $\widehat{\P}_N[.]=\P(.\mid S_N)$ 
 \end{lemma}
 \bpf
It is clear from the definition of $\xi^{\varphi(N)}_t$ that there exists a death process $(D_t, t\ge 0)$, which is independent of $(K_t, t\ge 0)$ conditionally upon $D_0=\varphi(N)-K_0$, and such that 
$$\xi^{\varphi(N)}_t=\tilde{K}_t+D_t, \ \forall t\ge 0,$$
where  
\begin{eqnarray*}
\tilde{K}_t=
\begin{cases}
K_t, & \;  \;    0\le t\le T_1;  \\
\bar{K_t}-1, &   \;  \;  t>T_1.\\
\end{cases}
\end{eqnarray*}
On the other hand, we have
$$\{\inf_{0\le t\le S_N}\xi^{\varphi(N)}_t>N\}\supset \{\inf_{0\le t\le S_N}D_t>N\}\supset\{D_{S_N}>N\}.$$
All we need to prove is that
$$\widehat{\P}_N(D_{S_N}\le N)\le \frac{2}{N^2}.$$
The process $(D_t, t\ge 0)$ is a jump Markov death process which takes values in the space $\{0, 1,\dots, \varphi(N)-K_0\}$. When in state $n$, $D_t$ jumps to $n-1$ at rate $\alpha n$ (recall that all crosses are kept in the present model). In other words the infinitesimal generator of $\{D_t, t\ge 0\}$ is given by 
$$Qf(n)=\alpha n[f(n-1)-f(n)].$$

 Let $f : \N\rightarrow \R$. The process $(M^{f}_t)_{t\ge 0}$ given by 
 \begin{equation}\label{Mar}
 M^f_t=f(D_t)-f(D_0)-\alpha \int_{0}^{t}D_s[f(D_s-1)-f(D_s)]ds 
 \end{equation}
 is a martingale. Applying \eqref{Mar} with the particular choice $f(n)=n$, there exists a martingale $(M^1_t)_{t\ge 0}$ such that $M^1_0=0$ and
 \begin{equation}\label{R-equation}
 D_t=D_0-\alpha \int_{0}^{t}D_s ds+M^1_t, \qquad t\ge 0.
 \end{equation}
 We note that $\{M^1_t, t\ge 0\}$ is a martingale under $\widehat{\P}_N[.]$. This is due to the fact that the Poisson process of crosses above $K_t$ is independent of $K_t$. We first deduce from \eqref{R-equation} that $\widehat{\E}_N(D_s)=D_0e^{-\alpha s}.$

 Using the fact that $D_t$ is a pure death process, we obtain the identity
 $$[M^1]_t=D_0-D_t,$$
which, together with \eqref{R-equation}, implies  
 \begin{equation*}\label{vq}
 <M^1>_t=\alpha \int_{0}^{t}D_s ds.
 \end{equation*} 
From \eqref{R-equation}, it is easy to deduce that (recall that $\varphi(N)=Ne^{\alpha S_N}(Ne^{\alpha S_N}+1)+K_0$)
 $$D_t=e^{-\alpha t}(\varphi(N)-K_0) +\int_{0}^{t}e^{-\alpha (t-s)}dM^1_s,$$
 which implies that 
 \begin{align*}
 \widehat{\P}_N(D_{S_N}\le N)&\le \widehat{\P}_N\Big( |\int_{0}^{S_N}e^{-\alpha (S_N-s)}dM^1_s |\ge N^2 e^{\alpha S_N}\Big)\\
 &= \widehat{\P}_N\Big( |\int_{0}^{S_N}e^{\alpha s}dM^1_s |\ge N^2 e^{2\alpha S_N}\Big)\\
  &\le\frac{1}{N^4e^{4\alpha S_N}}\int_{0}^{S_N}\alpha e^{2\alpha s}\widehat{\E}_N(D_s) ds\\
&\le \frac{2}{N^2}.
\end{align*}
The result is proved .   

\epf

\bigskip

From now on, we equip the probability space $(\Omega,\F,\P)$ with the filtration 
defined by $\F_t=\cap_{\eps>0}\mathring{\F}_{t+\eps}$, where  $\mathring{\F}_t=\sigma\{\eta_s(i),\ i\ge1,\ 0\le s\le t\}\vee{\mathcal N}$,
and $\mathcal N$ stands for the class of $\P$--null sets of $\F$. Any stopping time will be defined with respect to that filtration. 
\subsection{Exchangeability}
In this subsection, we will show that the $\Lambda$-lookdown model with selection preserves the exchangeability property, by an argument similar to that which we developed in \cite{bb.ep.bs}.

Let $S_n$ denote the group of permutations of the set $\{1, 2, \dots, n\}$. For all $\pi \in S_{n}$ and  $ a^{[n]}=(a_i)_{1\leq i  \leq n} \in \{0,1\}^{n}$, we define   the   vectors   
\begin{align*}
\pi^{-1}(a^{[n]})&=(a_{\pi^{-1}(1)},\ldots,a_{\pi^{-1}(n)}) =  (a^{\pi}_{i})_{1\leq i \leq n},\\
{\pi}(\xi^{[n]}_{t})&= (\xi_{t}(\pi(1)),\ldots ,\xi_{t}(\pi(n))).
\end{align*}

We   should   point  out   that   ${\pi}(\xi^{[n]}_{t})$ is a permutation of $(\xi_{t}(1),\ldots, \xi_{t}(n))$ and it is clear from the definitions that
\begin{equation}\label{relation-eta}
\{\pi(\xi^{[n]}_{t})=a^{[n]}\}=\{\xi^{[n]}_{t}=\pi^{-1}(a^{[n]})\}, \quad  \mbox{for any}\quad \pi  \in S_n.
\end{equation}

The main result of this subsection is 
\begin{theorem}\label{echange1}
If $(\eta_0(i))_{i\ge 1}$ are exchangeable random variables, then for all $t > 0$, $(\eta_t(i))_{i\ge 1}$ are exchangeable.
\end{theorem}
We first establish two lemmas, which treat repectively the case of resampling and of death events (we refer the reader to \eqref{measure} for the definition of the collection $\{t_i,\ i\ge1\}$).
\begin{lemma}\label{exchBirth}
For any finite stopping time $\tau$, any $\N$--valued $\F_{\tau}$--measurable random variable $n^*$,
if the random vector $\eta^{[n^*]}_{{\tau}}=(\eta_{{\tau}}(1),\ldots,\eta_{{\tau}}(n^*))$ is exchangeable, and
$T$ is the first time after $\tau$ of an arrow pointing to a level $\le n^*$ or a death at a level $\le n^*$,
then conditionally upon the fact that $T=t_{i_0}$, for some $i_0\ge 1$ and 
$\#(I_{t_{i_0}, p_{i_0}}\cap [n^*])=k$, where $k\ge 2$,
the random vector $\eta^{[n^*-1+k]}_{t_{i_0}}=\left(\eta_{t_{i_0}}(1),\ldots, \eta_{t_{i_0}}(n^*-2+k),\eta_{t_{i_0}}(n^*-1+k)\right)$ is exchangeable.
\end{lemma}
Note that $\eta^{[n^*-1+k]}_{t_{i_0}}$ is the list of the types of the individuals sitting
on levels $1,\ldots,n^*-1+k$ just after a birth event during which
one of the individuals sitting on a level between $1$ and $n^*$ has put $k-1$ children on levels up to $n^*$. 

\bpf 
For the sake of simplifying the notations, we condition upon $n^*=n, t_{i_0}=t$, $p_{t_{i_0}}=p$ and $\#(I_{t_{i_0}, p_{i_0}}\cap [n^*])=k$.
We start with some notation. 
\begin{align*}
A_{t}^{j_0,\dots,j_{k-1}}:=\{&\text{the $k$ levels selected by the point $(t, p)$ between levels $1$ and $n$ are}\\
& j_0, j_1, \dots, j_{k-1} \}. 
\end{align*}
We define
$$\widehat{\P}_{t,n} [.]=\P( . | t_{i_0}=t, n^*=n, \#(I_{t, p}\cap [n])=k). $$
Thanks  to  \eqref{relation-eta}, we   deduce  that, for $\pi\in S_{n-1+k}$, $a^{[n-1+k]}\in \{0, 1\}^{n-1+k}$,
\begin{equation}\label{relation-P}
\begin{split}
&\widehat{\P}_{t,n}(\pi(\eta^{[n-1+k]}_{t})=a^{[n-1+k]})\\
&=\sum_{1\leq j_0 <j_1<\dots<j_{k-1}\leq n}\!\!\widehat{\P}_{t,n}\left(\{\eta_{t}^{[n-1+k]}=(a^{\pi}_1,\ldots,a^\pi_{n-1+k})\}, A_{t}^{j_0,\dots,j_{k-1}}\right)
\end{split}
\end{equation}

On the event $A_{t}^{j_0,\dots,j_{k-1}}$, we have :
\begin{eqnarray*}
\eta_t(i)=
\begin{cases}
\eta_{t^{-}}(i),  & \mbox{if} \;  \;    1\leq i<j_1  \\
\eta_{t^{-}}(j_0),&   \mbox{if} \;  \;   i \in \{j_1, j_2, \dots, j_{k-1}\} \\
\eta_{t^{-}}(i-(\#\{\{j_1, j_2,\dots,j_{k-1}\}\cap [i]\})), &  \mbox{if} \;  \;    j_1<i\leq n-1+k, \;i\notin \{j_2, \dots, j_{k-1}\}.
\end{cases}
\end{eqnarray*}
This  implies  that   
$$A_{t}^{j_0,\dots,j_{k-1}}\cap\{\eta^{[n-1+k]}_{t}=(a^{\pi}_1,\ldots,a^\pi_{n-1+k})\}\subset\{ a^{\pi}_{j_0}=a^{\pi}_{j_1}=a^{\pi}_{j_2}=\dots=a^{\pi}_{j_{k-1}}\}.$$
 For $1< j_1<j_2<\dots<j_{k-1}\le n$, define the mapping $\rho_{j_1,j_2,\dots,j_{k-1}} : \{0, 1\}^{n+k-1}\longrightarrow \{0, 1\}^{n}$ by :
$$\rho_{j_1,j_2,\dots,j_{k-1}}(b_1, \ldots, b_{n-1+k})=(B_{j_1}, \dots, B_{j_{k-1}}), $$
where
\begin{align*}
B_{j_1}&=(b_{1},\ldots, b_{j_1-1}),\\
B_{j_m}&=(b_{j_{m-1}+1}, b_{j_{m-1}+2}, \dots, b_{j_{m} -1}) , \;\; 2\le m\le k-2\\
B_{j_{k-1}}&=(b_{j_{k-1}+1}, b_{j_{k-1}+2}, \dots, b_{n-1+k}).
\end{align*}
In other words, $\rho_{j_1,j_2,\dots,j_{k-1}}(z)$ is the vector $z$ from which the coordinates with indices $j_1,\dots,j_{k-1}$ have been suppressed. 
 The   right hand  side  of  \eqref{relation-P}  is  equal   to
\begin{equation*}
\sum_{1\leq j_0 <j_1<\dots<j_{k-1}\le n}{\bf 1}_{\{a^{\pi}_{j_0}=a^{\pi}_{j_1}\dots=a^{\pi}_{j_{k-1}}\}}\widehat{\P}_{t,n}\left(\{\eta^{[n]}_{t^{-}}=\rho_{j_1,j_2,\dots,j_{k-1}}(\pi^{-1}(a^{[n-1+k]}))\}, A_{t}^{j_0,\dots,j_{k-1}}\right).
\end{equation*}
It is easy to see that the events $(\eta^{[n]}_{t^{-}}=\rho_{j_1,j_2,\dots,j_{k-1}}(\pi^{-1}(a^{[n-1+k]})))$ and $A_{t}^{j_0,\dots,j_{k-1}}$ are independent. Thus
\begin{align*}
\widehat{\P}_{t,n}(\pi(\eta^{[n-1+k]}_{t})=a^{[n-1+k]})&=\sum_{1\leq j_0 <j_1<\dots<j_{k-1}\le n}{\bf 1}_{\{a^{\pi}_{j_0}=a^{\pi}_{j_1}\dots=a^{\pi}_{j_{k-1}}\}}\\
&\;\times\widehat{\P}_{t,n}\left(\eta^{[n]}_{t^{-}}=\rho_{j_1,j_2,\dots,j_{k-1}}(\pi^{-1}(a^{[n-1+k]}))\right)\widehat{\P}_{t,n}(A_{t}^{j_0,\dots,j_{k-1}})\\
&={n\choose k}^{-1}\sum_{1\leq j_0 <j_1<\dots<j_{k-1}\le n}{\bf 1}_{\{a^{\pi}_{j_0}=a^{\pi}_{j_1}\dots=a^{\pi}_{j_{k-1}}\}}\\
&\;\times\widehat{\P}_{t,n}\left(\eta^{[n]}_{t^{-}}=\rho_{j_1,j_2,\dots,j_{k-1}}(\pi^{-1}(a^{[n-1+k]}))\right).
\end{align*}
On the other hand, we have 
$$\#\{1\le j_0<\cdots<j_{k-1}\le n: a_{j_0}=\dots=a_{j_{k-1}}\}=\#\{1\le j_0<\cdots<j_{k-1}\le n : a^{\pi}_{j_0}\dots=a^{\pi}_{j_{k-1}}\}$$
Let $\ell_0<\ell_1<\cdots<\ell_{k-1}$ be the increasing reordering of the set $\{\pi(j_0),\pi(j_1), \cdots, \pi(j_{k-1})\}$. If $a_{j_0}=a_{j_1}=\cdots=a_{j_{k-1}}$, then we have $a^{\pi}_{\ell_0}=a^{\pi}_{\ell_1}\dots=a^{\pi}_{\ell_{k-1}}=a_{j_0}=a_{j_1}=\cdots=a_{j_{k-1}}$, and consequently $\rho_{j_1,j_2,\dots,j_{k-1}}(a^{[n-1+k]})$ and $\rho_{\ell_1,\ell_2,\dots,\ell_{k-1}}(\pi^{-1}(a^{[n-1+k]}))$ contain the same number of $0$'s and $1$'s. Since $\eta^{[n]}_{t^-}$ is exchangeable,
\begin{align*}
\widehat{\P}_{t,n}(\pi(\eta^{[n-1+k]}_{t})=a^{[n-1+k]})&={n\choose k}^{-1}\sum_{\gamma \in\{0, 1\}}\sum_{1\leq \ell_0 <\ell_1<\dots<\ell_{k-1}\le n}{\bf 1}_{\{a^{\pi}_{\ell_0}=a^{\pi}_{\ell_1}\dots=a^{\pi}_{\ell_{k-1}}=\gamma\}}\\
&\qquad \qquad \times\widehat{\P}_{t,n}\left(\eta^{[n]}_{t^{-}}=\rho_{\ell_1,\ell_2,\dots,\ell_{k-1}}(\pi^{-1}(a^{[n-1+k]}))\right)\\
&={n\choose k}^{-1}\sum_{\gamma \in\{0, 1\}}\sum_{1\leq j_0 <j_1<\dots<j_{k-1}\le n}{\bf 1}_{\{a_{j_0}=a_{j_1}\dots=a_{j_{k-1}}=\gamma\}}\\
&\qquad \qquad\times \widehat{\P}_{t,n}\left(\eta^{[n]}_{t^{-}}=\rho_{j_1,j_2,\dots,j_{k-1}}(a^{[n-1+k]})\right)\\
&=\widehat{\P}_{t,n}(\eta^{[n-1+k]}_{t}=a^{[n-1+k]}).
\end{align*}
The result follows.\epf
\begin{lemma}\label{exchDeath}
For any finite stopping time $\tau$, any $\N$--valued $\F_{\tau}$--measurable random variable $n^*$,
if the random vector $\eta^{[n^*]}_{{\tau}}=(\eta_{{\tau}}(1),\ldots,\eta_{{\tau}}(n^*))$ is exchangeable, and
$T$ is the first time after $\tau$ of an arrow pointing to a level $\le n^*$ or a death at a level $\le n^*$,
then conditionally upon the fact that $T$ is the time of a death, 
the random vector $\eta^{[n^*-1]}_{T}=\left(\eta_{T}(1),\ldots, \eta_{T}(n^*-1)\right)$ is exchangeable.\end{lemma}
\bpf 
To ease the notation we will condition upon $n^*=n$ and $T=t$. Let $\pi \in S_{n-1}$ 
and $a^{[n-1]}\in\{0,1\}^{n-1}$ be arbitrary. We consider the events :
$$B^i_t :=\{\text{the level of the dying individual at time $t$ is $i$} \}.$$ 
Let $\widehat{\P}_{t,n} [.]=\P(.|T=t,n^*=n)$.
We have
\begin{align*}
\widehat{\P}_{t,n}(\pi(\eta^{[n-1]}_{t})=a^{[n-1]})&=\sum_{1\leq i \leq n}\widehat{\P}_{t,n}\left(\eta^{[n-1]}_{t}=\pi^{-1}(a^{[n-1]}), B_{t}^{i}\right)\\
&=\sum_{1\leq i \leq n}\widehat{\P}_{t,n}\left(\eta_{t}(1)=a^{\pi}_{1},\ldots,\eta_{t}(n-1)=a^{\pi}_{n-1}, B_{t}^{i}\right).
\end{align*}
Define
$$ c^{\pi,n}_{i}=(a^{\pi}_{1},\ldots,a^{\pi}_{i-1}, 1, a^{\pi}_{i},\ldots, a^{\pi}_{n-1}),  \qquad c^{n}_{i}=(a_{1},\ldots,a_{i-1}, 1, a_{i},\ldots, a_{n-1}).$$
The last term   in  the    previous   relation    is equal  to 
\begin{align*}
\sum_{1\leq i \leq n}\widehat{\P}_{t,n}\left(\eta^{[n]}_{t^{-}}=c^{\pi,n}_i, B_{t}^{i}\right)
&=\sum_{1\leq i \leq n}\P\left(\eta^{[n]}_{t^{-}}=c^{\pi,n}_i\right)\widehat{\P}_{t,n}\left({B_{t}^{i}}\mid \eta^{[n]}_{t^{-}}=c^{\pi,n}_{i}\right)\\
&= \frac{1}{1+\sum_{j=1}^{n-1}a^{\pi}_{j} }\sum_{1\leq i \leq n}\P\left(\eta^{[n]}_{t^{-}}=c^{\pi,n}_{i}\right).
\end{align*}
Thanks  to   the exchangeability of $(\eta_{t^{-}}(1),\ldots,\eta_{t^{-}}(n))$, we have
\begin{equation*}
\widehat{\P}_{t,n}(\pi(\eta^{[n-1]}_{t})=a^{[n-1]})=\frac{1}{1+\sum_{j=1}^{n-1}a_{j}} \sum_{1\leq i \leq n}\P\left(\eta^{[n]}_{t^{-}}=c^{n}_{i}\right),
\end{equation*}
since  $\sum_{j=1}^{n-1}a^{\pi}_{j}=\sum_{j=1}^{n-1}a_{j}$ and $c_{i}^{\pi, n}$ is a permutation of $c_{i}^{\pi}$.  The result follows.
\epf

\bigskip

We can now proceed with the

\noindent{\sc Proof of Theorem \ref{echange1}}
For each $N\ge1$, let $\{V^N_t,\ t\ge0\}$ denote the $\N$--valued process which describes the position
at time $t$ of the individual sitting on level $N$ at time $0$, 
with the convention that, if that individual dies, we replace him by his neighbor below. The construction of our process $\{\eta_t(i), i\ge 1, t\ge 0\}$ in section 2.2 shows that $\inf_{t\ge0}V^N_t\to\infty$, as $N\to\infty$. 

It follows from Lemma \ref{exchBirth} and \ref{exchDeath} that for each $t>0$, $N\ge1$,
$(\eta_t(1),\ldots,\eta_t(V^N_t))$ is an exchangeable random vector.

Consequently, for any $t>0$, $n\ge1$, $\pi\in S_n$, $a^{[n]}\in\{0,1\}$,
$$|\P(\eta^{[n]}_t=a^{[n]})-\P(\eta^{[n]}_t=\pi^{-1}(a^{[n]}))|\le\P(V^N_t<n),$$
which goes to zero, as $N\to\infty$. The result follows.

\epf

\bigskip
For each $N\ge 1$ and $t\ge 0$, denote by $X^N_t$ the proportion of type $b$ individuals at time $t$ among the first $N$ individuals, i.e.
\begin{equation}
X^N_t=\frac{1}{N}\sum_{i=1}^{N}\eta_t(i)\label{XN}.
\end{equation} 
We are interested in the limit of $(X^N_t)_{t\ge 0}$ as $N$ tends to infinity. The following Corollary is a consequence of the well--known de Finetti's theorem (see e. g. \cite{Aldous}), which says that since they are exchangeable, the r.v.'s $\{\eta_t(i),\ i\ge1\}$ are i.i.d., conditionally upon their tail $\sigma$--field. 
\begin{corollary}\label{conv1}
For   each  $t\ge0$,
\begin{equation}\label{D-limit}
X_t =\lim_{N\rightarrow\infty}X_{t}^{N}\quad  \mbox{exists a.s.}
\end{equation}
\end{corollary}
\begin{remark}
Since the r.v.'s $\eta_t(i)$ take their values in $\{0,1\}$, their tail $\sigma$--field is exactly
$\sigma(X_t)$. 
 This fact will be used below.
\end{remark}

\section{Tightness and  Convergence to the $\Lambda$-W-F SDE with selection}\label{sectConv}
\subsection{Tightness of $\{X^N, N\ge 1, t\ge 0\}$}
In this part, we will prove the tightness of $(X^N)_{N\ge 1}$ in $D([0, \infty[)$, where for each $N\ge 1$ and $t\ge 0$, $X^N_t$ is defined by \eqref{XN}. 
For that sake, we shall write an integral equation for $X^N_t$.
We start with some notation.

For any $N, n, r, p$ such that $N\ge1, Nr\in \N, r\in ]0, 1], p\in [0, 1] $, we define \\
$Y( \cdot , N, p)$ to be the binomial distribution function with parameters $N$ and $p$;
$H(\cdot, N, n, r)$ the hypergeometric distribution function with parameters $(\!N-1, \!n-1, \!\frac{Nr-1}{N-1})$;
$\bar{H}(\!\cdot, \!N, n, r)$ the hypergeometric distribution function with parameters $(N-1, n-1, \frac{Nr}{N-1})$.
 For every $v, w \in [0, 1]$, let
\begin{align*}
F^N_p(v)&=\inf\{s; Y(s, N, p)\ge v\},\\
G_{N, n, r}(w)&=\inf\{s; H(s,N, n, r)\ge w\},\\
\bar{G}_{N, n, r}(w)&=\inf\{s; \bar{H}(s, N, n, r)\ge w\}.
\end{align*}
It follows that if $V, W$ are $\mathcal{U}([0, 1])$ r.v.'s, then the law of $F^N_p(V)$ is binomial with parameters $N, p$. $G_{N, n, r}(W)$(resp $\bar{G}_{N, n, r}(W)$) is hypergeometric with parameters $N-1, n-1, \frac{Nr-1}{N-1}$ (resp $N-1, n-1, \frac{Nr}{N-1}$).
Note that $F^N_p(\cdot)= Y^{-1}( \cdot , N, p), G_{N, n, r}(\cdot)=H^{-1}(\cdot, N, n, r)$ and $\bar{G}_{N, n, r}(\cdot)=\bar{H}^{-1}(\cdot , N, n, r)$. We recall that if $X$ is hypergeometric with parameters $(N, n, p)$ such that $Np\in \N$ and $p\in [0, 1]$, then
$$\E(X)=np\quad \text{and}\quad Var(X)=\frac{N-n}{N-1}np(1-p).$$
Now, for every $r, u, p, v, w \in[0, 1]$, let 
\begin{align}\label{psiN}
\psi^N(r, u, p, v, w)&=\frac{1}{N}{\bf 1}_{_{F^N_p(v)\ge 2}}\Big[{\bf 1}_{u\le r}\Big(F^N_p(v)-1-G_{N, F^N_p(v), r}(w)\Big)-{\bf 1}_{u> r}\bar{G}_{N, F^N_p(v), r}(w)\Big].
\end{align}
From the identity $r(n-1-\E[G_{N, n, r}(W)])=(1-r)\E[\bar{G}_{N, n, r}(W)]$, we deduce the
\begin{lemma}\label{TCPN}
For each $N\ge 1, r, p, v\in[0, 1]$ and $t\ge 0,$
$$\int_{]0, 1]^2 }\psi^N(r, u, p, v, w)dudw=0.$$
\end{lemma}
\bigskip

Using the definition of the model, one deduces that
\begin{align*}X^N_t=X^N_0&+\int_{ [0, t] \times]0, 1]^4 }\psi^N(X^N_{s^-}, u, p, v, w)M_0(ds, du, dp, dv, dw)\\
 &-\frac{1}{N}\int_{[0,t]\times[0,1]}{\bf1}_{u\le X^N_{s-}}{\bf1}_{\eta_{s-}(N+1)=0}M^N_1(ds,du)
\end{align*}
where $M_0$ and $M^N_1$ are two mutually independent Poisson point processes. $M_0$ is a Poisson point process on $\R_{+}\times [0, 1]\times[0,1]\times[0, 1]\times[0,1]$ with intensity measure $\mu(ds, du, dp, dv, dw)=dsdup^{-2}\Lambda(dp)dvdw$, $M^N_1$ is a Poisson point process on $\R_{+}\times[0, 1]$ with intensity measure $\alpha N\lambda(ds,du)=\alpha Ndsdu$. 
The reason why $X^N_t$ follows the above SDE is as follows.
 Births events happen according to the PPP $m$. With probability $X^N_{s^-}$, the individual which is copied (if at all) is of type 1. It is copied in a number which equals $(F^N_p(V)-1)^+$, where $F^N_p(V)$ follows the binomial law $(N,p)$. The increase in the number of 1's is that number, minus the number of ones which get pushed over level $N$, and that umber is the hypergeometric r.v. $G_{N,F^N_p(V),X^N_{s^-}}(W)$. In case 
the individual who is copied is a 0, the decrease in the number of ones is the hypergeometric r.v. $\bar{G}_{N,F^N_p(V),X^N_{s^-}}(W)$. Concerning the deaths, they 
happen according to a PPP with rate $\alpha NX^N_{s^-}$, and a death at time $s$ decreases the number of 1's by 1 iff $\eta_{s-}(N+1)=0$.

Now let 
\begin{equation}\label{M.Com}
\bar{M_0}=M_0-\mu, \quad \bar{M}^N_1=M^N_1-\alpha N\lambda.
\end{equation}
Using Lemma \ref{TCPN}, we have
\begin{equation}\label{ecritureXN}
\begin{split}
X^N_t&=X^N_0+\int_{ [0, t] \times]0, 1]^4 }\psi^N(X^N_{s^-}, u, p, v, w)\bar{M_0}(ds, du, dp, dv, dw)\\
&\qquad \ -\frac{1}{N}\int_{[0, t]\times[0, 1]}{\bf 1}_{u\le X^N_{s^-}}
{\bf 1}_{\eta_{s-}(N+1)=0}\bar{M}^N_1(ds, du)-\alpha\int_0^tX^N_s{\bf1}_{\eta_{s}(N+1)=0}ds.
\end{split}
\end{equation}

For each $N\ge 1, t\ge 0$, we define
\begin{align*}
\mathcal{M}^N_t&=\int_{ [0, t] \times]0, 1]^4 }\psi^N(X^N_{s^-}, u, p, v, w)\bar{M_0}(ds, du, dp, dv, dw)\\
\mathcal{N}^N_t&=\frac{1}{N}\int_{[0, t] \times]0, 1]}{\bf 1}_{u\le X^N_{s^-}}
{\bf 1}_{\eta_{s-}(N+1)=0}\bar{M}^N_1(ds, du)\\
V^N_t&=-\alpha\int_{0}^{t}X^N_s{\bf1}_{\eta_{s}(N+1)=0}ds.
\end{align*}
$\mathcal{M}^N_t$ and $\mathcal{N}^N_t$ are two orthogonal martingales. We have
\begin{equation*}\label{pro}
X^N_t=X^N_0+V^N_t+\mathcal{M}^N_t-\mathcal{N}^N_t.
\end{equation*}

$\forall N\ge 1, X^N_0 \in [0, 1] $, which implies that it is tight. Moreover, we have 
\begin{proposition}\label{Tightness}
The sequence $(X^N, N\ge 1)$ is tight in $D([0, \infty])$.
\end{proposition}

We first establish the lemma : 
\begin{lemma}\label{var1} For each $N\ge 1$ and $t\ge 0$,
\begin{align*}
\langle \mathcal{M}^N\rangle_t=\Lambda((0, 1)) \int_{0}^{t}X^N_s(1-X^N_s)ds\\
\langle\mathcal{N}^N\rangle_t=\frac{\alpha}{N} \int_{0}^{t}X^N_s{\bf1}_{\eta_{s}(N+1)=0}ds
\end{align*}
\end{lemma}
\bpf
Using the fact that $ \mathcal{M}^N$ and $ \mathcal{N}^N$ are pure-jump martingales, we deduce that 
$$\langle \mathcal{M}^N\rangle_t=\int_{ [0, t] \times]0, 1]^4 }(\psi^N(X^N_{s}, u, p, v, w))^2 dsdup^{-2}\Lambda(dp)dvdw.$$
Let 
\begin{align*}
\mathcal{A}^N(X^N_s, p)&=\int_{]0, 1]^3 }(\psi^N(X^N_{s}, u, p, v, w))^2 dudvdw\\
&=\frac{1}{N^2}\int_{]0, 1]^2}{\bf 1}_{_{F^N_p(v)\ge 2}}\Big[X^N_s\Big(F^N_p(v)-1-G_{N, F^N_p(v), X^N_{s^-}}(w)\Big)^2\\
&\quad \quad \quad \quad \quad \quad\quad\quad\quad+(1-X^N_s)(\bar{G}_{N, F^N_p(v), X^N_s}(w))^2\Big]dvdw.
\end{align*}
Tedious but standard calculations yield
\begin{align*}
\int_{ [0, 1] }\Big[r\Big( F^N_p(v)-1-G_{N, F^N_p(v), r}(w)\Big)^2dw+(1-r)\Big(\bar{G}_{N, F^N_p(v), r}(w)\Big)^2\Big]dw\\
=\frac{N}{N-1}r(1-r)F^N_p(v)(F^N_p(v)-1),
\end{align*}
for every $v, r \in [0, 1]^2$.
Consequently
\begin{align*}
\mathcal{A}^N(X^N_s, p)&=\frac{X^N_s(1-X^N_s)}{N(N-1)}\int_{[0, 1]}{\bf 1}_{F^N_p(v)\ge 2}F^N_p(v)(F^N_p(v)-1)dv\\
&=p^2X^N_s(1-X^N_s).
\end{align*}
We deduce that
\begin{align*}
\langle \mathcal{M}^N\rangle_t &=\int_{ [0, t]\times[0, 1]}\mathcal{A}^N(X^N_s, p)ds p^{-2}\Lambda(dp)\\
&=\Lambda((0, 1))\int_{ [0, t]}X^N_s(1-X^N_s)ds.
\end{align*}
Similarly, we have
\begin{equation*}\label{Mart2}
\begin{split}
\langle \mathcal{N}^N\rangle_t&=\frac{\alpha}{N}\int_{[0, t]\times[0, 1]}{\bf 1}_{u\le X^N_{s}}{\bf 1}_{\eta_s(N+1)=0}dsdu\\
&=\frac{\alpha}{N}\int_{[0, t]}X^N_s{\bf 1}_{\eta_s(N+1)=0}ds.
\end{split}
\end{equation*}
The lemma has been established.
\epf
\bigskip

We can now proceed with the

\noindent{\sc Proof of Proposition \ref{Tightness}}
We have
$$
X^N_t=X^N_0+V^N_t+\mathcal{M}^N_t-\mathcal{N}^N_t$$ and
$$\langle \mathcal{M}^N-\mathcal{N}^N\rangle_t=\langle \mathcal{M}^N\rangle_t+\langle\mathcal{N}^N\rangle_t.$$
Moreover, from Lemma \ref{var1}
\begin{equation*}\label{major}
\left|\frac{dV^N_t}{dt}\right|\le\alpha,\quad 0\le \frac{d \langle \mathcal{M}^N\rangle_t}{dt}\le\frac{\Lambda((0,1))}{4},\quad
0\le\frac{d \langle \mathcal{N}^N\rangle_t}{dt}\le\frac{\alpha}{N}. 
\end{equation*} 
Aldous' tightness criterion (see Aldous \cite{A.T}) is an easy consequence of 
those estimates. 
\epf

\bigskip
Now, from Proposition \ref{Tightness} and \eqref{D-limit}, it is not hard to show there exists a process $X\in D([0, \infty))$, such that for all $t\ge 0,$
\begin{equation}\label{X-limit}
X^N_t \rightarrow X_t\ \ a.s, 
\end{equation}
and
$$X^N\Rightarrow X\; \mbox{weakly in}\; D([0, \infty)).$$

\subsection{Convergence to the $\Lambda$-Wright-Fisher SDE with selection}
Our goal is to get a representation of the process $(X_t)_{t\ge0}$ defined in \eqref{X-limit} as the unique weak solution to the stochastic differential equation \eqref{eqXt}.

Let $(\Omega, {\mathcal F}, \;  \P)$ be a fixed probability space, on which the above Poisson measures are defined, which is  equipped with the filtration described at the end of section 2.2. Recall the Poisson point measure $M=\sum_{i=1}^{\infty}\delta_{t_i, u_i, p_i}$ defined in the Introduction, and  
for every $u\in ]0, 1[$ and $r\in [0, 1]$, we introduce the elementary function
$$\Psi(u, r)={\bf 1}_{u\le r}-r.$$
We rewrite equation \eqref{eqXt} as
\begin{equation}\label{LWF}
X_t=x- \alpha \int_{0}^{t} X_s(1-X_s)ds+\int_{[0, t]\times]0, 1[^2}p\Psi(u, X_{s^-})\bar{M}(ds, du, dp), \;   t>0, 0<x<1,
\end{equation}
which we call the $\Lambda$-Wright-Fisher SDE with selection. Without loss of generality, we shall assume that $\alpha>0$, which means that $X_t$ represents the proportion of non-advantageous alleles.

The proof of the following identity is standard and left to the reader.  

\begin{lemma}\label{conv2} For each $r\in[0, 1]$,
\begin{align*}
\int_{[0, 1]^4}&\left(\psi^N(r, u, p, v, w)-p\Psi(u, r)\right)^2dup^{-2}\Lambda(dp)dvdw\\
&=2r(1-r)\Big[\frac{N}{N-1}\int_{[0, 1]^2}(1-up)^{N-1}du\Lambda(dp)-\frac{\Lambda([0, 1])}{N-1}\Big].
\end{align*}
\end{lemma}

Let us now prove the main result of this section.
\begin{theorem}\label{main-result}
Suppose that $X^N_0\rightarrow x$ a.s., as $N\rightarrow \infty$. Then the $[0,1]-$valued process $\{X_t ,t \ge 0\}$ defined by \eqref{X-limit} is the unique solution to the  $\Lambda$-W-F SDE with selection \eqref{LWF}.
\end{theorem}
\bpf
Strong uniqueness of the solution to \eqref{LWF} follows from Theorem 4.1 in \cite{DL}. We now prove that $(X_t)_{t\ge 0}$ defined by \eqref{X-limit} is a solution to the $\Lambda$-Wright-Fisher \eqref{LWF}.

We know that  $X^N_t\rightarrow X_t$ a.s. for all $t\ge0$ and that $X^N\Rightarrow X$ weakly in $D([0,+\infty))$ as $N\to\infty$. 
Recall the decomposition
\begin{equation}\label{decomp}
X^N_t=X^N_0+\mathcal{M}^N_t-\mathcal{N}^N_t+V^N_t.
\end{equation}
It follows from Lemma \ref{var1} that $\mathcal{N}^N_t\to0$ in probability, as $N\to\infty$.
We next show that
\begin{equation}
\begin{split}
\mathcal{M}^N_t\rightarrow 
\int_{[0, t]\times]0, 1[\times]0, 1[}p\Psi(X_{t^-}, u)\bar{M}(ds, du, dp)\ \text{in probability},\quad as \;N\rightarrow \infty\label{lim2},
\end{split}
\end{equation}
where $\overline{M}=\int_{[0,1]^2}\bar{M_0}(., ., ., dv, dw)$.
For each $N\ge 1$ and $t\ge 0$, let
$$h^N(t)=\int_{[0, t]\times[0, 1]^4}\left(\psi^N(X^N_{s^-}, u, p, v, w)-p\Psi(X_{s^-}, u)\right)\bar{M_0}(ds, du, dp, dv, dw),$$
where $\bar{M_0}$ is defined by \eqref{M.Com}.
 $\{h^N(t),\ t\ge0\}$ is a martingale, and
\begin{align*}
\langle h^N\rangle_t=\int_{[0, t]\times[0, 1]^4}\left(\psi^N(X^N_{s}, u, p, v, w)-p\Psi(X_{s}, u)\right)^2dsdup^{-2}\Lambda(dp)dvdw.
\end{align*}
We have 
\begin{align*}
\langle h^N\rangle_t &\le 2\int_{[0, t]\times[0, 1]^2}\left(p\Psi(X^N_{s},u)-p\Psi( X_{s},u)\right)^2dsdup^{-2}\Lambda(dp)\\
&+ 2t\sup_{0\le s\le t}\int_{[0, 1]^4}\left(\psi^N(X^N_{s}, u, p, v, w)-p\Psi( X^N_s,u)\right)^2dup^{-2}\Lambda(dp)dvdw\\
&\le 2\int_{[0, t]\times[0, 1]^2}\left(p\Psi(X^N_{s},u)-p\Psi(X_{s},u)\right)^2dsdup^{-2}\Lambda(dp)\\
&+ 2t\sup_{0\le r\le 1}\int_{[0, 1]^4}\left(\psi^N(r, u, p, v, w)-p\Psi( r,u)\right)^2dup^{-2}\Lambda(dp)dvdw.
\end{align*}
Using the fact that $X^N_s\rightarrow X_s$ a.s., it is not hard to show by the dominated convergence theorem that as $N\rightarrow \infty$,
\begin{equation}\label{c1}
\int_{[0, t]\times[0, 1]^2} \left(\Psi(X^N_{s},u)-\Psi(X_{s},u)\right)^2dsdu\Lambda(dp)\rightarrow 0 \quad \mbox{a.s}.
\end{equation}
Now from lemma \ref{conv2} , it is easy to show that as $N\rightarrow \infty$,
\begin{equation}\label{c2}
\sup_{0\le r\le 1}\int_{[0, 1]^4}\left(\psi^N(r, u, p, v, w)-p\Psi( r,u)\right)^2dup^{-2}\Lambda(dp)dvdw\rightarrow 0
\end{equation}
Combining \eqref{c1} and \eqref{c2}, we deduce that
$$\forall t\ge 0,\quad\;\langle h^N\rangle_t \rightarrow 0\quad a.s, \quad as \; N\rightarrow \infty.   $$
On the other hand, from the bound $|\psi(r,u)|\le1$ and Lemma \ref{conv2}, we deduce that 
$$\langle h^N\rangle_t\le Ct \Lambda([0, 1]),\; \forall N\ge 2.$$ 
Hence from the dominated convergence theorem 
$$\lim_{N\rightarrow\infty}\E [h^N(t)^2]=0 \quad\;\forall t\ge 0$$
i.e.
\begin{align*}
\mathcal{M}^N_t=\int_{[0, t]\times[0, 1]^4}\psi^N(X^N_{s^-}, u, p, v, w)\bar{M_0}(ds, du, dp, dv, dw)&\xrightarrow{L^2} \int_{[0, t]\times[0, 1]^4}p\Psi(X_{s^-}, u)\bar{M}(ds, du, dp)
\end{align*}
as $N\rightarrow \infty$, in particular
$$\mathcal{M}^N_t\rightarrow\int_{[0, t]\times[0, 1]^2}p\Psi(X_{s^-}, u)\bar{M}(ds, du, dp)\; \text{in probability}\; ,as \;N\rightarrow \infty\ .$$
\eqref{lim2} is established. 

From \eqref{decomp}, we deduce that
$$\frac{1}{N}\sum_{k=1}^N V^k_t=\frac{1}{N}\sum_{k=1}^NX^k_t-\frac{1}{N}\sum_{k=1}^NX^k_0-\frac{1}{N}\sum_{k=1}^N\mathcal{M}^k_t+\frac{1}{N}\sum_{k=1}^N\mathcal{N}^k_t.$$
It follows from the above arguments and our assumption on the initial condition that for all $t\ge0$, as $N\to\infty$, the right--hand side converges in probability towards
$$X_t-x-\int_{[0, t]\times[0, 1]^4}p\Psi(X_{s^-}, u)\bar{M}(ds, du, dp).$$
But
\begin{align*}
\frac{1}{N}\sum_{k=1}^N V^k_t&=-\alpha\int_0^t\frac{1}{N}\sum_{k=1}^NX^k_s{\bf1}_{\eta_s(k+1)=0}ds\\
&=-\alpha\int_0^t\frac{1}{N}\sum_{k=1}^N\left(X^k_s-X_s\right){\bf1}_{\eta_s(k+1)=0}ds
-\alpha\int_0^tX_s\left(\frac{1}{N}\sum_{k=1}^N{\bf1}_{\eta_s(k+1)=0}\right)ds\\
&\to0-\alpha\int_0^tX_s(1-X_s)ds
\end{align*}
a.s., as $N\to\infty$. The result follows clearly from the above facts.
\epf

\begin{remark}
Our proof establishes in fact that for all $t>0$, as $N\to\infty$, 
$$\int_0^tX^N_s{\bf1}_{\eta_s(N+1)=0}ds\to\int_0^tX_s(1-X_s)ds$$
in probability. This does not mean that ${\bf1}_{\eta_s(N+1)=0}$ converges, but it seems intuitively clear that for any $0\le r<t$,
$$\int_r^t{\bf1}_{\eta_s(N+1)=0}ds\to\int_r^t(1-X_s)ds.$$
However, that convergence is not really easy to establish.
\end{remark}

\begin{remark} Suppose we know a priori that $(X_t)_{t\ge0}$ defined by \eqref{X-limit} is a Markov process. Then we can prove that $(X_t)_{t\ge0}$ is a solution to the $\Lambda$-Wright-Fisher SDE \eqref{LWF}
as follows. Let us look backwards from time $t$ to time $0$. For each $0 \le s \le t$, we denote by $Z^{n,t}_s$ the highest level occupied by the ancestors at time $s$ of the $n$ first individuals at time $t$. We know that conditionally upon $X_t$, the $\{\eta_t(i), i\ge 1\}$ are i.i.d Bernoulli with parameter $X_t$. Consequently, for any $n\ge 1$,
$$X^n_t=\P(\eta_t(1)=\dots=\eta_t(n)=1\mid X_t),$$
this implies that
\begin{align*}
\E_x[X^n_t]&=\E_x[\P(\eta_t(1)=\dots=\eta_t(n)=1\mid X_t)]\\
&=\P_x(\eta_t(1)=\dots=\eta_t(n)=1)\\
&=\P_x(\text{the $1\dots Z^{n, t}_0$ individuals at time $0$ are all $b$})\\
&=\E_n[x^{Z^{n, t}_0}].
\end{align*}  
It is plain that the conditional law of $Z^{n, t}_0,$ given that $(\eta_t(1)=\dots=\eta_t(n)=1)$ equals the conditional law of $R_t$, given that $R_0=n$. Consequently, for each $n\ge 1$
$$\E[X^n_t\mid X_0=x]=\E[Y^n_t\mid Y_0=x],$$
where $(Y_t)_{t\ge 0}$ is a solution to \eqref{LWF}. But for all $t > 0, r \in [0, 1]$, the conditional law of $X_t$, given that $X_0 = x$ is determined by its moments, since $X_t$ is a bounded r. v. So 
$(X_t)_{t\ge0}$ and $(Y_t)_{t\ge0}$ have the same transition densities, that is $\{X_t, t\ge 0\}$ is the unique weak solution to \eqref{LWF}.
\end{remark}

\subsection{An alternative proof of uniqueness}

Uniqueness in law could also by proved as in \cite{B.LG} (where the case $\alpha=0$ is treated) by a duality argument, which we now sketch .

Recall the notation $\Psi(u, y)={\bf 1}_{u\le y}-y$. For every $y\in[0, 1]$ and every function $g : [0, 1]\rightarrow \R$ of class $\mathcal{C}^2$, we set
$$\mathcal{L}g(y)=\int_{[0, 1]\times[0, 1]} \Big[g(y+p\Psi(u, y))-g(y)-p\Psi(u, y)g^{\prime}(y)\Big]p^{-2}\Lambda(dp)du-\alpha g'(y)(1-y)y.$$ 
A solution $(Y_t)_{t\ge 0}$ of \eqref{LWF} is a Markov process with generator $\mathcal{L}$.
Hence for every $g:[0,1]\rightarrow\R$ of class $\mathcal{C}^2$, the process
\begin{equation*}
g(Y_t)-\int_{0}^{t}ds\mathcal{L}g(Y_s), \quad \quad t\ge 0 
\end{equation*}
is a martingale.


It is plain that for $ g(z)=z^n$
\begin{equation}\label{gen}
\mathcal{L}g(z)=\sum_{k=2}^{n}{n\choose k}\lambda_{n, k}(z^{n-k+1}-z^n)+\alpha n(z^{n+1}-z^n).
\end{equation}

Let $\{{ R}_t, t\geq 0\}$  be   a $\N$-valued jump Markov process  which, when in state $k$, jumps to 
\begin{enumerate}
\item $k-\ell+1$ at rate  ${k\choose \ell} \lambda_{k, \ell} , \quad 2\le \ell \le k$; 
\item $k+1$  at rate $\alpha k$, $\alpha>0$.
\end{enumerate}

In other words,
the infinitesimal generator of $\{{R}_t, t\geq 0\}$ is given by:
$$\mathcal{L}^{*}f(k)=\sum_{\ell=2}^{k} {k\choose \ell} \lambda_{k, \ell}[f(k-\ell+1)-f(k)]+\alpha k[f(k+1)-f(k)].$$
For every $z\in [0, 1]$ and every $r\in\N$, we set
\begin{equation}\label{r-duality}
P(z, r )=z^{r}.
\end{equation}
Viewing $P(z, r )$ as a function of $r$, we have
\begin{align*}
\mathcal{L}^{*}P(z, r )&=\sum_{k=2}^{r} {r\choose k} \lambda_{r, k}[z^{r-k+1}-z^{r}]+\alpha r[z^{r+1}-z^{r}]
\end{align*}
On the other hand, viewing $P(z, r )$ as a function of $z$ we can easily evaluate $\mathcal{L}P(z, r )$ from formula \eqref{gen}, and we deduce that 
\begin{equation}
\mathcal{L}P(z, r )=\mathcal{L}^{*}P(z, r ).\label{gen1=gen2}
\end{equation}
Now suppose that $(Y_t)_{t\ge 0}$ is a solution to \eqref{LWF}, and let $R_0=n$. By a standard argument (see Section 4.4 in \cite{E.K}) we deduce from \eqref{gen1=gen2}  that
$$\E[P(Y_t, R_0)]=\E[P(Y_0, R_t)],$$ i.e
\begin{equation*}\label{duality}
\E [Y_t^{n}|Y_0=x]=\E[x^{R_t}|R_0=n].
\end{equation*}
Since this is true for each $n\ge 1$ and $Y_t$ take values in the compact set [0, 1], this is enough to identify the conditional law of $Y_t$, given that $Y_0=x$, for all $0\le x\le 1$. 
Since $(Y_t)_{t\ge0}$ is a homogeneous Markov process, this implies that the law of $(Y_t)_{t\ge0}$ is uniquely determined.

\section{Fixation and non-fixation in the $\Lambda$-W-F SDE}
\subsection{The CDI property of the $\Lambda$-coalescent}
In this subsection, we recall a remarkable property of the $\Lambda$-coalescent $(\Pi_t)_{t\ge 0}$  defined in the introduction. For each $n\ge 1$, let $\#\Pi^{[n]}_t$ denote the number of blocks in the partition $\Pi^{[n]}_t$ ($\Pi^{[n]}_t$ is the restriction of $\Pi_t$ to $[n]$ ). Then let $T_n=\inf\{t\ge 0 : \#\Pi^{[n]}_t=1\}$. As stated in (31) of \cite{Pitman}, we have 
$$0=T_1<T_2\le T_3\le\dots \uparrow T_{\infty}\le \infty.$$ 
We say the $\Lambda$-coalescent comes down from infinity ($\Lambda \in$ \textbf{CDI}) if $\P(\# \Pi_t<\infty)=1$ for all $t>0$, and we say it stays infinite if $\P(\#\Pi_t=\infty)=1$ for all $t>0$.  The coalescent comes down from infinity if and only if $T_{\infty}<\infty$ a.s. We will show that this is equivalent to fixation. Kingman showed that the $\delta_{0}$-coalescent comes down from infinity.

A necessary and sufficient condition for a $\Lambda$-coalescent to come down from
infinity was given by Schweinsberg \cite{J.Sch1}. Define
$$
\phi(n)=\sum_{k=2}^{n}(k-1){n\choose k}\lambda_{n, k},
$$
and 
$$\nu(dp)=p^{-2}\Lambda(dp).$$
It is not hard to deduce from the binomial formula that
$$\phi(n)=\int_{0}^{1}[np-1+(1-p)^n]\nu(dp).$$
Schweinsberg's result \cite{J.Sch1} says that the $\Lambda$-coalescent comes down from infinity if and only if 
\begin{equation}\label{CDI}
\sum_{n=2}^{\infty}\frac{1}{\phi(n)}<\infty.
\end{equation}  We shall see below that the convergence of this series is also necessary and sufficient for fixation in finite time. Using the fact that the function $f_n(p)=(1-p)^n-1$ is decreasing for any fixed $n$, we have
\begin{equation*}\label{Major.gamma}
\int_{0}^{1}(np-1)\nu(dp)\le\phi(n)\le n\int_{0}^{1}p\nu(dp), \qquad \forall n\ge 1.
\end{equation*}
The last assertion together with \eqref{CDI}, implies that if $\int_{0}^{1}p\nu(dp)<\infty$ then the $\Lambda$-coalescent stays infinite. This result has been proved by Pitman (see lemma 25 in \cite{Pitman}).
 


Theorem \ref{main-result} shows that $(X_t)_{t\ge0}$ is a bounded supermartingale. Indeed, if $(X_t)_{t\ge 0}$ is a solution to \eqref{LWF} , then for all $0\le t\le s$, 
\begin{align*}
\E(X_t\mid {\mathcal F}_s)&\le x- \alpha \int_{0}^{s} X_r(1-X_r)dr+\E\Big[\int_{[0, t]\times]0, 1[\times]0, 1[}p\Psi(u, X_{s^-})\bar{M}(ds, du, dp)\mid {\mathcal F}_s\Big]\\
&=X_s.
\end{align*} 
Consequently the following limit exists a.s
\begin{equation}\label{X-infty}
X_{\infty}=\lim_{t \rightarrow \infty}X_t \in \{0, 1\}.
\end{equation}
Indeed, 0 and 1 are the only possible limit values. 
\subsection{Fixation and non-fixation in the $\Lambda$-W-F SDE }
We assume that the initial proportion $x$ of type $B$ individuals satisfies $0<x<1$.
In this section, we prove that  fixation happens in finite time iff the condition \eqref{CDI} is satisfied.  
Before establishing the main result of this section, we collect some results which will be required for its proof.
\begin{lemma}\label{lim}
$$ \frac{\phi(n)}{n}\uparrow\int_0^1p\nu(dp)\;\; \mbox{as} \;\; n \uparrow \infty,$$
where $\nu(dp)=p^{-2}\Lambda(dp)$.

 \end{lemma}
 \bpf 
 \begin{align*}
 \phi(n)&=\int_{0}^{1}\left[np-1+(1-p)^n\right]\nu(dp)\\
 &=\int_{0}^{1}\left[\frac{n}{p}\left(1-\int_{0}^{1}(1-up)^{n-1}du\right)\right]\Lambda(dp).
 \end{align*}
On the last line, we have made use of the identity
 $$(1-p)^n-1=\int_{0}^{1}-np(1-up)^{n-1}du.$$
For each $p \in ]0, 1]$, let $$f^n(p)=\frac{1}{p}\left(1-\int_{0}^{1}(1-up)^{n-1}du\right).$$ We have,
 \begin{align*}
 n^{-1} \phi(n)=\int_{0}^{1}f^n(p)\Lambda(dp).
 \end{align*}
 The result follows from the monotone convergence theorem.
 
\epf

We now deduce that
\begin{lemma}\label{var-phi}
The function $\phi$ increases, and 
$$\sum_{n=2}^{\infty}\frac{1}{\phi(n)}<\infty\quad \Rightarrow\quad \sum_{n=2}^{\infty}\frac{1}{\phi(n)-\alpha n}<\infty.$$ 
\end{lemma}
\bpf
We have 
\begin{align*}
\phi(n+1)-\phi(n)&=\int_{0}^{1}[p+(1-p)^{n+1}-(1-p)^n]\nu(dp)\\
&=\int_{0}^{1}p(1-(1-p)^n)\nu(dp)\\
&\ge 0.
\end{align*}
Which implies the first claim. Now, we already know that if $\sum_{n=2}^{\infty}\frac{1}{\phi(n)}<\infty$, then $\int_0^1p\nu(dp)=\infty$. Thus, the second assertion is a consequence of the last lemma and the following relation
$$\sum_{n=2}^{\infty}\frac{1}{\phi(n)-\alpha n}=\sum_{n=2}^{\infty}\frac{\alpha}{\phi(n)(n^{-1}\phi(n)-\alpha)}+\sum_{n=2}^{\infty}\frac{1}{\phi(n)}.
$$
The lemma is proved.

\epf
\bigskip

For each $t\ge 0$, we define again
$$K_t=\inf\{i\ge 1 : \eta_t(i)=0\}.$$ and $$T_1=\inf\{t\ge 0 : K_t=1\}.$$ 
We have the following 
\begin{theorem}\label{Fixation}
If $\Lambda \in$ \textbf{CDI}, then one of the two types ( $b$ or $B$) fixates  in finite time, i.e. 
$$\exists\; \zeta <\infty\; a.s : X_{\zeta}=X_{\infty} \in \{0, 1\}  $$
If $\Lambda \notin$ \textbf{CDI}, then
$$\forall t\ge0,\; 0<X_t<1 \;a.s.$$
\end{theorem} 
\bpf
The proof has been inspired by \cite{Bertoin} (see Section 4 ). 


\textsc{Step 1} : Suppose that $\Lambda \in$ \textbf{CDI}.  We consider two cases.  \\
\textsc{Case 1} : $K_0=1$.\\ In this case, the allele $B$ fixates in the population. Indeed, the individual at level 1 never dies and he cannot be pushed to an upper level. Let 
$$\zeta=\inf\{t> 0 : \eta_t(i)=0, \forall i\ge 1\} .$$
$\zeta$ is the time of fixation of allele $B$. We are going to show that $\zeta<\infty\ a.s$. 

 We couple our original population process with the following $\N$--valued process
   $\{Y_t,\ t\ge0\}$, which describes the growth of a population which we denote 
   ``the $B$--population'', and whose dynamics  we now describe.
   $Y_0=1$, at time zero the $B$--population consists of a unique individual who occupies site 1, while all other sites $k\ge2$ are empty. We follow the same realizations of the Poisson point process $m$ on $\R_+\times[0,1]$ (see \eqref{measure})  and of the sets $I_{t,p}$ as presented in the Introduction. 
   
   At each time $t$ corresponding to an atom $(t,p)$ of the Poisson point process $m$, 
   we associate the set $I_{t,p}$. We put a cross at time $t$ on
   all levels $i\in I_{t,p}$, except the lowest one. If there is at least one cross on the interval $[2,Y_{t-}+1]$, we modify the population as follows
   (otherwise we do nothing).  All individuals sitting at time $t^-$ below the lowest cross don't move. All others are displaced 
   upwards in such a way that all sites with a cross become free, and the respective orders of the individuals remain unchanged. Finally, individuals are added 
   on all sites with a cross which  lie below or immediately above an occupied site.  Clearly, as long as the growing number of individuals of the $B$--population remains below any given value $k$, the number of atoms of the Poisson process $m$ which modify the size of the population on any given finite time interval remains finite, and each jump in the population size is finite.
    However, we will now show that as a consequence of the CDI property of the associated coalescent process, the jumps of $Y_t$ accumulate in such a way that 
   $Y_t=+\infty$, for some finite (random) $t$. Since it is plain that $Y_t$ is less than the total number of type $B$ individuals in the population, this will show that 
   $\zeta<\infty$ a.s.
   
   Indeed, looking backward in time, starting from any $t>0$ the process which describes 
   the genealogy of the ``$B$--population'' is the Lambda--coalescent. More precisely, as a time--reversal of our $B$--population
   process, it is the Lambda--coalescent starting from the random value $Y_t$, and conditioned upon the fact that all the partitions 
   have coalesced into one single partition by time 0. 
   
   This claim is justified as follows.  
   Let $\{U_s=Y_{t-s},\ 0\le s\le t\}$. At each time $s$ of a point $(s,p)$ of the PPP $m$ where $\sharp\left(I_{s,p}\cap[1,U_{s-}]\right)\ge2$,
   all lineages of the set $I_{s,p}\cap[1,U_{s-}]$ coalesce. Would we describe the evolution of $\{U_s,\ s\ge0\}$ using copies of $m$ and the $I_{s,p}$'s which would be 
   independent of those used to describe the growth of $Y_\cdot$, then $U_\cdot$, starting from $U_0=Y_t=N$, would be an instance of the $N$--$\Lambda$--coalescent. Here and below  we make a slight abuse of terminology, calling $\Lambda$--coalescent the process which describes the number of blocks in a
   $\Lambda$--coalescent. 
   
  For each $N\ge2$, we define
  $$\xi_N=\inf\{t>0,\ Y_t\ge N\},$$
  and by $\theta_N$ the time taken by the $N$--$\Lambda$--coalescent to reach 1. It follows from an obvious coupling that $N\to\theta_N$ is increasing.
  In fact we shall only use the fact that $N\to\E\theta_N$ is increasing. Since $Y_{\xi_N-}<N$, it is plain that
   \begin{equation}\label{eq:ineqfin}
  \E\xi_N\le\E\xi_2+\E\theta_N,
  \end{equation}
  and moreover the law of $\xi_2$ is exponential with parameter $\int_0^1p^2\nu(dp)$.
  Let us admit for a moment the
  \begin{lemma}\label{le:bound}
  For any $N>1$, $$\E\theta_N\le\sum_{k=2}^\infty\phi(k)^{-1},$$
and this bound is finite since $\Lambda\in\textbf{CDI}$. 
  \end{lemma}
Since $\zeta\le\lim_{N\to\infty}\xi_N$, it follows from \eqref{eq:ineqfin} and Lemma \ref{le:bound} that
$\E\zeta<\infty$.

In order to conclude Case 1 of the first step of the proof of our Theorem, let us proceed with the

\noindent{\sc Proof of Lemma \ref{le:bound}}
The Markov process which describes the number of ancestors in a $\Lambda$--coalescent jumps from $n$
   to $n-\ell+1$ ($2\le\ell\le n$) at rate $\begin{pmatrix}n\\ \ell\end{pmatrix} \lambda_{n,\ell}$ . In other words, its infinitesimal generator $Q$ is given by
   $$Qf(n)=\sum_{\ell=2}^n  \begin{pmatrix}n\\ \ell\end{pmatrix} \lambda_{n,\ell}[f(n-\ell+1)-f(n)].$$
Let us define for each $n\ge 1$
$$f(n)=\sum_{k=n+1}^{\infty}\frac{1}{\phi(k)}.$$
We have for $2\le \ell \le n$
$$f(n-\ell+1)-f(n)=\sum_{k=n-\ell +2}^{n}\frac{1}{\phi(k)}.$$
Recall Lemma \ref{var-phi}. Since $1/\phi$ is decreasing, we have for $2\le \ell \le n$,  
$$f(n-\ell+1)-f(n)\ge (\ell -1)\frac{1}{\phi(n)},$$
and therefore
$$
Qf(n)\ge \frac{1}{\phi(n)}\sum_{\ell=2}^{n} {n\choose \ell}(\ell -1) \lambda_{n, \ell}=1.
$$
Using the fact that the process 
$$f(U_t)-f(U_0)-\int_{0}^{t}Qf(U_s)ds, \qquad t\ge 0$$
is a martingale, we obtain 
\begin{align*}
\E(\theta_N)&\le \E\left(\int_{0}^{\theta_N}Qf(U_s)ds\right)\\
&\le f(1)
\end{align*}
\epf

\bigskip

{\sc Case 2} : $K_0>1$. 

If $T_1<\infty$ then type $B$ fixates in finite time. Indeed, wait until $T_1$ which is a stopping time at
which the Markov process $\{\eta_t(i),i\ge 1\}_{t\ge 0}$ starts afresh, and then use the 
argument from Case 1.  

We suppose now that $T_1=\infty$, which implies that $K_t\rightarrow \infty$ as $t\rightarrow \infty$, as already noted in section 2.2.2. In other words, if $T_1=\infty$, then the allele $B$ does not fixate in the population. Let 
$$n_0=\inf\{n\ge 1 : \phi(n)-\alpha n\ge 1\}.$$
Such an $n_0$ exists because since $\Lambda \in$ \textbf{CDI}, $\int_{0}^{1}p\nu(dp)=+\infty$, hence by Lemma \ref{lim}, we have $\lim_{n\rightarrow \infty}n^{-1}\phi(n)=+\infty$. 

We define a ``$b$--population'' $\{Y_t,\ t\ge0\}$, which again starts from a unique ancestor sitting on level 1.
The novelty is that now each individual dies at rate $\alpha$. It then may happen that the ``$b$--population'' gets empty. In that case, we immediately start afresh with a new unique ancestor sitting at level 1.  The fact that eventually the ``$b$--population'' grows
and become larger than any $N$ is a consequence of the fact that $K_t\rightarrow \infty$ as $t\rightarrow \infty$.

Note that the process describing the number of ancestors of the present individuals in that population  is now a jump--Markov process with generator
$Q_\alpha$ given by
$$Q_\alpha f(n)=\sum_{\ell=2}^{n} {n\choose \ell} \lambda_{n, \ell}[f(n-\ell+1)-f(n)]+
\alpha n (f(n+1)-f(n)),$$
conditioned upon hitting 1 before time $t$.

Let  $N>n_0$ denote a fixed integer, $\xi_N$ the time taken by the ``$b$--population'' to reach the value $N$, i.e.
   $$\xi_N=\inf\{t>0,\ Y_t\ge N\},$$
   and by $\theta^{n_0}_N$ the time taken by the process with generator $Q_\alpha$ to come down below $n_0$, starting from $N$.
   Similarly as in \eqref{eq:ineqfin}, we have 
  \begin{equation}\label{eq:ineq2} 
  \E\xi_N\le\E\xi_{n_0}+\E\theta^{n_0}_N.
  \end{equation}
   In order to show that  the allele $b$ fixates in finite time, it remains to establish the
   \begin{lemma}
   There exists a constant $C<\infty$ such that 
   $$\E\theta^{n_0}_N\le C,$$
   for all $N>n_0$.
   \end{lemma}
   \bpf
 For each $n\ge 1$, we define
$$f(n)=\sum_{k=n+1}^{\infty}\frac{1}{(\phi(k)-\alpha k)\vee 1}.$$
By Lemma \ref{var-phi}, for each $n\ge 2$, $f(n)$ is finite. We have for $2\le \ell \le n$
$$f(n-\ell+1)-f(n)=\sum_{k=n-\ell +2}^{n}\frac{1}{(\phi(k)-\alpha k)\vee 1}.$$
Since $k\rightarrow 1/(\phi(k)-\alpha k)\vee 1$ is decreasing, we obtain  
$$f(n-\ell+1)-f(n)\ge (\ell -1)\frac{1}{(\phi(n)-\alpha n)\vee 1},$$
and therefore
\begin{align*}
Qf(n)&\ge \frac{1}{(\phi(n)-\alpha n)\vee 1}\sum_{\ell=2}^{n} {n\choose \ell}(\ell -1) \lambda_{n, \ell}-\frac{\alpha n }{(\phi(n+1)-\alpha (n+1))\vee 1}\\
&=\frac{\phi(n)}{(\phi(n)-\alpha n)\vee 1}-\frac{\alpha n }{(\phi(n+1)-\alpha (n+1))\vee 1}\\
&\ge \frac{\phi(n)}{(\phi(n)-\alpha n)\vee 1}-\frac{\alpha n }{(\phi(n)-\alpha n)\vee 1}
\end{align*}
hence $Qf(n)\ge 1, $ for each $n\ge n_0$. Since the process 
$$f(U_t)-f(U_0)-\int_{0}^{t}Qf(U_s)ds, \qquad t\ge 0$$
is a martingale and $U_t$ remains bounded while $0\le t\le \theta^{n_0}_N$, 
\begin{align*}
\E(\theta^{n_0}_N)&\le \E\left(\int_{0}^{\theta^{n_0}_N}Qf(U_s)ds\right)\\
&=f(U_{\theta^{n_0}_N})-f(U_0)\\
&\le f(1)
\end{align*}
\epf

\textsc{Step 2} : Suppose $\Lambda \notin$ \textbf{CDI}, that is the $\Lambda$-coalescent does not come down from infinity. We have 
\begin{equation}\label{NCDI}
\sum_{n=2}^{\infty}\frac{1}{\phi(n)}=+\infty.
\end{equation}
We claim that $(K_t,t \ge 0)$ does not reach $\infty$ in finite time. The contrary would imply that $\exists T<\infty$ such that $K_{T}=\infty$ a.s., so the number of ancestors at tiome 0 of the infinite population at time $T$ in the $\Lambda$-lookdown model would be finite, which contradicts the fact that $\Lambda \notin$ \textbf{CDI }. 
Hence $K_t<\infty$ a.s. This implies that $X_t<1$, for all $t\ge 0$. Indeed if $X_t=1$, for some $t>0$, by applying de Finetti's Theorem, we deduce that $\eta_t(i)=1, \forall i\ge 1$, which contradicts the fact that $K_t<\infty$. It remains to show that $X_t>0$ for all $t\ge 0$.

For any $m\ge 1, t>0$, we define the event
$$ A^m_t=\{\text{The $m$ first individuals of type $b$ at time 0 are dead at time $t$} \}$$
We have
$$\P(A^m_t)=(1-e^{-\alpha t})^m,$$
and then 
$$\P(\cap_m A^m_t)=0 \quad \ \forall t> 0. $$
From this, we deduce that $\exists i\ge 1$ such that $\eta_t(i)=1$. 
The same argument used for the proof of $X_t<1$ now shows that $X_t>0$, for all $t\ge 0$.
\epf

\subsection{The law of $X_{\infty}$}
Let $x$ be the proportion of type $b$ individuals at time $0$, where $0<x<1$. As the individual at level 1 cannot be pushed to an upper level, we have
$$\{\eta_0(1)=0\}\subset \{X_{\infty}=0\}, \; \mbox{hence}\;  \P(X_{\infty}=0)\ge 1-x.$$
If $\alpha =0$, $(X_t)_{t\ge0}$ is a bounded martingale, so 
$$\P(X_{\infty}=1)=\E(X_{\infty})=\E(X_0)=x.$$ 
If $\alpha >0$, by using \eqref{LWF} together with \eqref{X-infty}, we deduce that
$$\P(X_{\infty}=1)=\E X_{\infty}<x.$$
In this subsection we want to describe those cases where can we decide whether $\P(X_{\infty}=1)>0$ or $\P(X_{\infty}=1)=0$. We first prove 

\begin{proposition}
If $\Lambda \in \textbf{CDI}$, then 
$$\P(X_{\infty}=1)>0.$$
\end{proposition}
\bpf
Since $\Lambda \in \textbf{CDI},$ if all individuals at time $0$ would be of type $b$, there would be a (random) level $J$ such that the individual sitting on level $J$ at time $0$ reaches $+\infty$ in finite time. Now $\P(X_{\infty}=1)>0$ follows from the fact that $\P(K_0>J)>0$, where $K_0$ denotes the lowest level occupied by a type $B$ individual at time $0$.   
\epf
\bigskip

In the case $\Lambda \not\in \textbf{CDI}$, since selection has infinite time to act, one may wonder whether or not $\P(X_\infty=1)=0$.
Some partial results have been obtained in that direction in Bah \cite{BB}, but since then the question has been completely settled by 
Foucart \cite{CF} and Griffiths \cite{RG}, who prove
\begin{theorem}\label{FouGri}
Suppose that $0<x<1$, and let $$\alpha^\ast:=\int_0^1\log\left(\frac{1}{1-p}\right)\nu(dp).$$
\begin{enumerate}
\item If $\alpha< \alpha^\ast$, then $0<\P(X_\infty=0|X_0=x)<1$.
\item If $\alpha\ge\alpha^\ast$, then $X_\infty=0$ a.s.
\end{enumerate}
\end{theorem}
Needless to say, if $\alpha^\ast=+\infty$, which is in particular the case when $\Lambda \in \textbf{CDI}$, we are in the first case.
Note that \cite{CF} settles the two cases $\alpha<\alpha^\ast$ and $\alpha>\alpha^\ast$, while \cite{RG} treats the case $\alpha=\alpha^\ast$ as well,
assuming $\alpha^\ast<\infty$
in the first case.
We refer to \cite{CF} and \cite{RG} for references to earlier partial results on this problem in the biological literature.

\subsection{The fixation line, special case of the Bolthausen--Sznitman coalescent}\label{ss:hn}

The aim of this section is to connect our model and results with the recent work of H\'enard \cite{Henard}, 
and to compute the law of $X_\infty$ and the speed at which either type invades the whole population, in the case
of the Bolthausen--Sznitman coalescent.


H\'enard's definition of the fixation line is as follows. Consider the levels of the offsprings at time $t>0$ of the individual sitting at time 0 at level 1. This constitutes a subset of $\N$, the connected component containing 1 of which is of the form $\{1,\ldots,L_t\}$. This defines the fixation line $L_t$. In our case (in contradiction with H\'enard's situation), there may be no such connected component containing 1, if $\eta_0(1)=1$ and $\eta_t(1)=0$ for some $t>0$, in which case we define $L_t$ to be 0. H\'enard's fixation line is an increasing process. Our is increasing if the individual sitting
on level 1 at time 0 is of type $B$ (i.e. is a 0), but this is not the case if that individual is of type $b$ (i. e. is a 1).

We are only interested in this second case, which is the only one where conditionally upon the value of $\eta_0(1)$, $X_\infty$ is random.
However, we will not necessarily assume that $L_0=1$. We prefer to define the fixation line as follows.

For all $t\ge0$, let
$$L_t=\max\{k\ge1;\ \eta_1(t)=\eta_2(t)=\cdots=\eta_k(t)=1\},$$
and this defines also $L_0$. Equivalently, $L_t=K_t-1$, where $K_t$ is the lowest level occupied at time $t$ by an individual of type 0, see
the discussion in subsection \ref{Llms}. $L_t$ is clearly a $\Z_+$--valued continuous time Markov process.

$L_t$ does not evolve as discussed
in \cite{Henard}, since those individuals sitting on levels $\{1,2,\ldots, L_0\}$, as well as their offsprings, are type $b$ individuals, who die at rate $\alpha$, each death inducing a jump of $L_t$ 
of size $-1$.
The process $\{L_t,\ t\ge0\}$ is a $\Z_+$--valued Markov process, whose jump rates are given by
$$\Gamma_{i,j}=\begin{cases}
\begin{pmatrix} j\\ j-i+1\end{pmatrix}\int_0^1x^{j-i-1}(1-x)^i\Lambda(dx),&\text{if $1\le i<j<\infty$};\\
\alpha i,&\text{if $j=i-1$};
\end{cases}
$$
whenever $i\ge1$, and the process is absorbed at 0. Indeed, $\Gamma_{i,j}$ is the rate at which $K_t$ jumps from $i+1$ to $j+1$.
As was shown in subsection \ref{Llms}, either $L_t\to\infty$, as $t\to\infty$, in which case $X_\infty=1$, or else $L_t$ hits zero in finite time, in which case $X_\infty=0$. 
In the first case, $L_t$ explodes in finite time iff $\Lambda\in \textbf{CDI}$. In the case where  $\Lambda\not\in \textbf{CDI}$,
it is of interest to describe the speed at which $L_t\to\infty$, whenever this happens. This is done in the case without selection (and it applies in our situation to the case where $K_0=1$)
in \cite{Henard}, in the situation $\Lambda(dx)=dx$, i.e. the case of the Bolthausen--Sznitman coalescent. We will show that the same result applies in our case, i.e. the 
slow--down due to the death essentially does not modify that speed. 

Recall that the Bolthausen--Sznitman coalescent  belongs to the family of the Beta$(2-\alpha,\alpha)$ ($0<\alpha<2$) coalescents, it corresponds to the case $\alpha=1$. Note that the Beta$(2-\alpha,\alpha)$ coalescent comes down from infinity iff $1<\alpha<2$. The Bolthausen--Sznitman coalescent is the border case. 
One may expect that in this model, on the event $\{L_t\to\infty\}$,  $L_t\to\infty$ very fast, as $t\to\infty$.

Before going to that, let us compute explicitly the law of $X_\infty$, in that case of the Bolthausen--Sznitman coalescent. 
The possibility of that computation is due to the remark that in this particular case (and only in that one), the process $L_t$
is a continuous time branching process. Indeed
in the case $\Lambda(dx)=dx$, we have
$$\Gamma_{i,i+j}=\begin{cases}
\frac{i}{j(j+1)},&\text{if $j\ge1$};\\
\alpha i,&\text{if $j=-1$};
\end{cases}
$$
This means that $L_t$ is a Markov continuous time branching process, with life time exponential with parameter $1+\alpha$, and family size distribution 
$\{p_j,\ j=0,2,3,\ldots\}$ given by $$p_0=\frac{\alpha}{1+\alpha},\quad  p_j=\frac{1}{j(j-1)(1+\alpha)}.$$
Note that the generating function of that probability distribution is given by
$$h(s)=\frac{s+\alpha}{1+\alpha}+\frac{1-s}{1+\alpha}\log(1-s).$$
We have
\begin{proposition}\label{prop:extinct-proba}
Conditionally upon $L_0=k$ ($k\ge1$),
$$\P(L_t=0)=\left[1-\exp\{-\alpha(1-e^{-t})\}\right]^k,\quad \P(\lim_{t\to\infty}L_t=0)=\left[1-e^{-\alpha}\right]^k.$$
\end{proposition}
\bpf
It suffices to consider the case $L_0=1$, which we now do. In that case, it follows from general results on continuous time branching processes,
see e.g. chapter V in \cite{TH}, that the collection of generating functions $f_t(s)=\E[s^{L_t}]$  satisfies the ODE
\begin{align*}
\partial_tf_t(s)&=\Phi(f_t(s)),\\
f_0(s)&=s,
\end{align*}
where $\Phi(z)=(1+\alpha)(h(z)-z)=(1-z)[\alpha+\log(1-z)]$ is the so--called infinitesimal generating function. 
It is not too hard to check that the solution of that ODE is
$$f_t(s)=1-\exp\left[\alpha(e^{-t}-1)+e^{-t}\log(1-s)\right].$$
Hence $$\P(L_t=0)=f_t(0)=1-\exp\left[\alpha(e^{-t}-1)\right],$$
from which the result follows. \epf

We can now conclude
\begin{corollary}\label{lawXinfty}
Again in the case $\Lambda(dx)=dx$,  
$$\P(X_\infty=0|X_0=x)=\frac{1-x}{1-x(1-e^{-\alpha})}.$$
\end{corollary}
\bpf  Recall that $L_0=K_0-1=0$ iff level 1 is occupied by a type $B$ individual at time 0, and that
at time 0 individuals  placed at levels $1,2,\ldots$ are choosen in an i.i.d. manner, each one being of type $b$ (i.e. 1) with probability $x$, and of type $B$ 
(i.e. 0) with probability $1-x$.
We have 
\begin{align*}
\P(X_\infty=0|X_0=x)&=\sum_{k=0}^\infty (1-e^{-\alpha})^k\P(L_0=k)\\
&=(1-x)\sum_{k=0}^\infty[x(1-e^{-\alpha})]^k\\
&=\frac{1-x}{1-x(1-e^{-\alpha})},
\end{align*}
were we have used Proposition \ref{prop:extinct-proba} for the first equality.
The result follows. \epf

Note that in the case $\Lambda(dx)=dx$, Theorem \ref{FouGri} tells us that $0<\P(X_\infty=0)<1$ for all $\alpha>0$ since $\alpha^\ast=+\infty$, which is consistent with 
the last result. Note also that \cite{RG} gives, for a general $\Lambda$ coalescent, an expression for the above quantity in terms of the sum of an infinite series.
It does not seem easy to deduce our result from that formula.
\begin{remark}
The proportion of advantageous alleles is $Y_t=1-X_t$. Our formula says (here ``BS'' refers to the Bolthausen--Sznitman coalescent)
$$\P_{\text{BS}}(Y_\infty=1|Y_0=y)=\frac{y}{y+(1-y)e^{-\alpha}}.$$
If we replace the Bolthausen--Sznitman by Kingman's coalescent, it is well--known (see e. g. \cite{RG}) that
the formula reads
$$\P_{\text{K}}(Y_\infty=1|Y_0=y)=\frac{1-e^{-2\alpha y}}{1-e^{-2\alpha}}.$$
We note that these two formulae coincide, and are equal to $(1+e^{-\alpha})^{-1}$, in the case $y=1/2$.
The following comparison holds : for all $\alpha>0$, $\P_{\text{BS}}(Y_\infty=1|Y_0=y)>\P_{\text{K}}(Y_\infty=1|Y_0=y)$ if $0<y<1/2$, while 
$\P_{\text{BS}}(Y_\infty=1|Y_0=y)<\P_{\text{K}}(Y_\infty=1|Y_0=y)$ if $1/2<y<1$. Indeed the difference 
$\P_{\text{BS}}(Y_\infty=1|Y_0=y)-\P_{\text{K}}(Y_\infty=1|Y_0=y)$ has the same sign as
$$\Phi(y)=e^{-2\alpha y}(e^{-\alpha}+y(1-e^{-\alpha}))+(e^{-\alpha}-e^{-2\alpha}) y-e^{-\alpha}.$$
Now $\Phi(0)=\Phi(1/2)=\Phi(1)=0$, $\Phi'(0)>0$, $\Phi'(1/2)<0$ and $\Phi'(1)>0$ for all $\alpha>0$, while $\Phi''(y)$ vanishes at the unique point
$$0<y_\alpha=\frac{1-e^{-\alpha}-\alpha e^{-\alpha}}{\alpha(1-e^{-\alpha})}<1.$$
\end{remark}

We now establish
\begin{theorem}\label{vitconv}
In the case $\Lambda(dx)=dx$, if $L_0=1$, then conditionally upon $L_t\to\infty$ as $t\to\infty$, 
$$e^{-t}\log L_t\to\mathbf{e}\quad\text{a.s.},$$
where $\mathbf{e}$ is a standard exponential r.v. 
\end{theorem}
\bpf
Recalling the infinitesimal generating function $\Phi$ specified in Proposition \ref{prop:extinct-proba},  it is not hard to see that the function 
$$\frac{1}{\Phi(1-x)}-\frac{1}{x\log x}$$ is integrable near zero (one way to see that is to make the change of variable $y=1/x$, and note that the resulting
integral, say from 2 to $\infty$, converges, by comparison with a Bertrand series). Hence condition (3) of Theorem 3 from Grey \cite{Grey} is satisfied, which 
implies the stated convergence, but is remains to specify the law of $\mathbf{e}$. 

We follow the strategy of proof of Proposition 3.8 in \cite{Henard}. For each $t>0$, $s\mapsto f_t(s)$ is a bijection from $[0,1]$ onto $[1-\exp\{-\alpha(1-e^{-t})\},1]$.
Its inverse reads 
$$g_t(s)=1-\exp\left[\alpha(e^t-1)+e^t\log(1-s)\right].$$ It
 is a bijection  from $[1-\exp\{-\alpha(1-e^{-t})\},1]$ onto $[0,1]$.
For each $1-\exp\{-\alpha\}\le s\le1$, $0<g_t(s)\le1$, and the process $\{ g_t(s)^{L_t},\ t\ge0\}$ is Markov and has constant expectation. Indeed
  $$\E\left[g_t(s)^{L_t}\right]=f_t(g_t(s))=s.$$
  Hence it is a $[0,1]$--valued martingale, which converges a.s. as $t\to\infty$ to a r.v. $V(s)$.
  Moreover, by dominated convergence and explicit computation, for any $\beta>0$,
$$\E[V(s)^\beta]=\lim_{t\to\infty}\E\left[g_t(s)^{\beta L_t}\right]=\lim_{t\to\infty}f_t[g_t(s)^\beta]=s.$$
This implies that $V(s)$ takes values in $\{0,1\}$, and $\P(V(s)=1)=\E(V(s))= s$.

Let us now define the r.v.
$$U=\inf\{1-e^{-\alpha}<s\le1,\ V(s)=1\}$$
It is plain that $\{U\le s\}=\{V(s)=1\}$, hence $\P(U\le s)=s$, for $s\in(1-e^{-\alpha},1]$. On the other hand, since $g(1-e^{-\alpha})=1-e^{-\alpha}$, we have that
$\{V(1-e^{-\alpha})=1\}=\{L(t)\to0\}$, and we see that the law of $U$ has a Dirac measure of mass $1-e^{-\alpha}$ at $1-e^{-\alpha}$,
and has density $1$ on the interval $(1-e^{-\alpha},1)$.

For $0\le s\le e^{-\alpha}$, we have that $\log[g_t(1-s)]\simeq - (\rho s)^{e^t}$ as $t\to\infty$, where $\rho=e^\alpha$.
If $s<1-U$, then 
\begin{align*}
g_t(1-s)^{L_t}&\to1, \text{ hence}\\
L_t\log[g_t(1-s)]&\to0, \text{ or equivalently}\\
(\rho s)^{e^t}L_t&\to0, \text{ as }t\to\infty,
\end{align*}
while if $s>1-U$,
\begin{align*}
g_t(1-s)^{L_t}&\to0, \text{ hence}\\
L_t\log[g_t(1-s)]&\to-\infty, \text{ or equivalently}\\
(\rho s)^{e^t}L_t&\to\infty, \text{ as }t\to\infty.
\end{align*}
Let $\Theta:=1-U$. For any $\eps>0$,
$$[\rho(\Theta-\eps)]^{e^t}L_t\to0,\ \text{ and } [\rho(\Theta+\eps)]^{e^t}L_t\to+\infty,\text{ as }t\to\infty.$$
Taking again the logarithm, we deduce that as $t\to\infty$,
\begin{align*}
\alpha e^t + \log(\Theta-\eps)e^t+\log(L_t)&\to-\infty,\\
\alpha e^t + \log(\Theta+\eps)e^t+\log(L_t)&\to+\infty.
\end{align*}
Consequently
$$-\log(\Theta+\eps)-\alpha\le\liminf_{t\to\infty}e^{-t}\log(L_t)\le\limsup_{t\to\infty}e^{-t}\log(L_t)\le-\log(\Theta-\eps)-\alpha.$$
This being true for any $\eps>0$, we have proved that, as $t\to\infty$,
$$e^{-t}\log(L_t)\to-\log(\Theta)-\alpha.$$
We now define $\mathbf{e}:=-\log(\Theta)-\alpha$. Conditionally upon $L_t\to\infty$ as $t\to\infty$, the law of $U$ is uniform on $[1-e^{-\alpha},1]$,
hence the law of $\Theta$ is uniform on $[0,e^{-\alpha}]$. Then for $r>0$,
$$\P(\mathbf{e}>r|L_t\to\infty)=\P(\Theta<\exp[-(\alpha+r)]|L_t\to\infty)=\exp[\alpha]\exp[-(\alpha+r)]=e^{-r}.$$
\epf

Let $\beta^n$ the time taken by the fixation line $L_t$, starting from $L_0=1$, to exceed the value $n$. As noted in \cite{Henard}, a consequence of 
Theorem \ref{vitconv} is that
$$\beta^n-\log\log(n)\to-\log(\mathbf{e})\ \text{ a.s. as }n\to\infty.$$
 In the situation treated in \cite{Henard}, $\beta^n$ has the same law as $\tau^n_1$, the time taken 
by the $n$--$\Lambda$ coalescent to hit the value 1, i.e. the time taken for $n$ individuals to find their most recent common ancestor.

In our case, the $\Lambda$--coalescent must be replaced by the $\Lambda$ -- Ancestral Selection Graph. Indeed, since in the forward time direction individuals die, in the backward time direction we have birth of lineages.

The $n$--$\Lambda$--ASG is defined as follows. Starting from $n$ lineages, the lineages coalesce according  to the $\Lambda$--coalescent, while 
new lineages are born according to the following rule. While there are $k\ge2$ active lineages, 
a new lineage is born at rate $\alpha k$, this lineage being placed on a level chosen uniformly among the levels 
$\{1,2,\ldots,k+1\}$. If the level $\ell\le k$ is chosen, the lineages located on levels $\ell,\ell+1,\ldots,k$ just before the birth event get pushed one level up. 
We refer to \cite{KN}  and \cite{NK} for the description of the ASG, where the coalescent is Kingman's coalescent. We note that here we consider only type $b$ individuals, type $B$
individuals occupying possibly some of the higher levels. 

Define $\tau^n_1$ to be the time for the $n$--$\Lambda$--ASG to find a common ancestor, i.e. the time for the number of lineages to reduce to 1. 
It follows from Theorem \ref{vitconv} that, in the case $\Lambda(dx)=dx$, as $n$ gets large, the decrease of the number of lineages due to the coalescence events is much faster than the creation of new lineages, hence $\tau^n_1<\infty$ a.s.
\cite{Henard} shows that in the case $\alpha=0$, the law of $\tau^n_1$ coincides with that of $\beta^n$, the time taken by the fixation line starting from 1 to reach a value greater than or equal to $n$. This is no longer true in the case $\alpha>0$, since the process of the number of lineages in the $n$--$\Lambda$--ASG is no longer decreasing. Here $\beta^n$ has rather the law of 
the time elapsed between the last time when there are at least $n$ lineages in the $n$--$\Lambda$--ASG, and the time when there is one lineage.
However for large $n$ this does not make a real difference, as follows
from the following result.
\begin{lemma}
Fix an arbitrary $h>0$. On the event that $L_t\to\infty$ as $t\to\infty$, for $n$ large enough, $L_{\tau^n_1+s}>n$, for all $s\ge h$.
\end{lemma}
\bpf Choose $\varepsilon>0$ small enough such that 
$$\frac{\mathbf{e}-\varepsilon}{\mathbf{e}+\varepsilon}>e^{-h}.$$
It follows from Theorem \ref{vitconv} that there exists $t_\varepsilon$ such that for any $t\ge t_\varepsilon$,
$$\mathbf{e}-\varepsilon\le e^{-t}\log L_t\le \mathbf{e}+\varepsilon,$$
and $n_\varepsilon$ such that whenever $n\ge n_\varepsilon$, 
$$\log\log n-\log(\mathbf{e}+\varepsilon)\le \tau^n_1\le\log\log n-\log(\mathbf{e}-\varepsilon).$$
Choose $n\ge n_\varepsilon$ such that moreover $\tau^n_1\ge t_\varepsilon$. Consequently
$$e^{\tau^n_1}\ge\frac{\log n}{\mathbf{e}+\varepsilon},$$
and whenever $s\ge h$,
\begin{align*}
L_{\tau^n_1+s}&\ge\exp\left(e^{\tau^n_1+s}(\mathbf{e}-\varepsilon)\right)\\
&\ge\exp\left(e^{\tau^n_1}e^s(\mathbf{e}-\varepsilon)\right)\\
&>\exp\left(e^{\tau^n_1} (\mathbf{e}+\varepsilon)\right)\\
&\ge n.
\end{align*}
\nobreak\epf

Consequently, the time elapsed between the first visit of a level above $n$ by $L_t$, and the last visit below $n$ after that time (if any) 
tends to zero in probability, as $n\to\infty$.
As a result, we can conclude as in \cite{Henard}
\begin{proposition}
Suppose we are again in the case $\Lambda(dx)=dx$, and define $\tau^n_1$ as above. Then, as $n\to\infty$,
$$\tau^n_1-\log\log n\Rightarrow -\log \mathbf{e}.$$
\end{proposition}

\begin{remark} We expect that our look--down construction, and the duality with the $\Lambda$--ASG can produce new results
beyond the case of the Bolthausen--Sznitman coalescent, at least in the case of the Beta--coalescents, in particular concerning the law 
of the number of blocks implied in the last coalescence in the Beta$(2-\alpha,\alpha)$--ASG, and the expectation of the depth of the Beta$(2-\alpha,\alpha)$--ASG
in case $1<\alpha<2$.
\end{remark} 

\section{Kingman and $\Lambda$--coalescent}
In this last section we suppose that the measure $\Lambda$ is general (i.e $\Lambda(\{0\})>0$). This implies that $\nu$ is infinite. Note that we could have $\Lambda((0,1))=0$, but this case corresponds to ``pure Kingman'', which is already well understood, see in particular \cite{bb.ep.bs}. So we assume again that \eqref{condition} is satisfied. We will show that the proportion $X_t$ of type $b$ individuals at time $t$ in the population of infinite size is a solution to the stochastic differential equation with selection
\begin{align}\label{X_t-general}
\begin{split}
X_t = x &- \alpha \int_{0}^{t} X_s(1-X_s)ds\\
&+\int_{0}^{t}\sqrt {cX_s(1-X_s)}dB_s\\
&+\int_{[0, t]\times]0, 1[\times]0, 1[}p({\bf 1}_{u\le X_{s^-}}-X_{s^-})\bar{M}(ds, du, dp),
\end{split}
\end{align}
where $c=\Lambda(\{0\}), \bar{M}$ is the compensated measure $M$ defined in section 3.2, and $B$ is a standard Brownian motion. Let $\{W(ds, du)\}$ be a white noise on $(0, \infty)\times (0, 1]$ based on the Lebesgue measure $dsdu$. We remark that if $(X_t)_{t\ge0}$ satisfies \eqref{X_t-general}, then $X_t$ is a solution in law of the following stochastic differential equation
\begin{align*}
\begin{split}
X_t = x&- \alpha \int_{0}^{t} X_s(1-X_s)ds\\
&+\sqrt {c}\int_{[0, t]\times]0, 1[}({\bf 1}_{u\le X_s}-X_s)
W(ds, du)\\
&+\int_{[0, t]\times]0, 1[^2}p({\bf 1}_{u\le X_{s^-}}-X_{s^-})\bar{M}(ds, du, dp).
\end{split}
\end{align*}   

We first define the model. Recall the process $\{\eta_t(i), i\ge 1, t\ge 0\}$ defined in the introduction.
The evolution of the population is the same as that described in the case $\Lambda(\{0\})=0$ except that we superimpose single births, which are described as follows 
\begin{itemize}
\item[\hspace{12pt}] 
For any $1\le i<j$, arrows are placed from $i$ to $j$ according to a rate $\Lambda(\{0\})$ Poisson process, independently of the other pairs $i'<j'$. Suppose there is an arrow from $i$ to $j$ at time $t$. Then a descendent (of the same type) of the individual sitting on level $i$ at time $t^-$ occupies the level $j$ at time $t$, while for any $k\ge j$, the individual occupying the level $k$ at time $t^-$ is shifted to level $k+1$ at time $t$. In other words, $\eta_t(k)=\eta_{t^-}(k)$ for $k<j$, $\eta_t(j)=\eta_{t^-}(i)$, $\eta_t(k)=\eta_{t^-}(k-1)$ for $k>j$. 
\end{itemize}
By coupling our model with the simplest lookdown model with selection defined in \cite{bb.ep.bs}, it is not hard to show that for $N$ large enough, the individual sitting on level $2N$ at time 0 never visits a level below $N$, that is the evolution within the box $(t, i)\in[0, \infty)\times \{1, 2, \dots, N\}$ is not altered by removing all crosses above $2N$. The process $\{\eta_t(i), i\ge 1, t\ge 0\}$ is well-defined.  


For each $N\ge 1$ and $t\ge 0$, denote by $X^N_t$ the proportion of type $b$ individuals at time $t$ among the first $N$ individuals, i.e.
\begin{equation}
X^N_t=\frac{1}{N}\sum_{i=1}^{N}\eta_t(i)\label{XNN}.
\end{equation}
Combining the arguments in \cite{bb.ep.bs} and section 2.3 (see above), it is easy to show if $(\eta_0(i))_{i\ge 1}$ are exchangeable random variables, then for all $t > 0$, $(\eta_t(i))_{i\ge 1}$ are exchangeable. An application of de Finetti's theorem, yields that
\begin{equation}\label{NX-limit}
X_t =\lim_{N\rightarrow\infty}X_{t}^{N}\quad  \mbox{exists a.s.}
\end{equation}

Using the definition of the model, it is easy to see that ($\psi^N$ was defined by \eqref{psiN})
\begin{align*}
X^N_t=X^N_0+\mathcal{K}^N_t&+\int_{ [0, t] \times]0, 1]^4 }\psi^N(X^N_{s^-}, u, p, v, w)\bar{M}(ds, du, dp, dv, dw)\\
&-\frac{1}{N}\int_{[0, t]\times[0, 1]}{\bf 1}_{u\le X^N_{s^-}}{\bf 1}_{\eta_{s-}(N+1)=0}M^N_1(ds, du),
\end{align*}
where $\mathcal{K}^N_t$ is a martingale of jump size $\pm\frac{1}{N}$. We have 
\begin{lemma}\label{crochM}
$$\langle \mathcal{K}^N\rangle_t=\int_{0}^{t}\varphi^N(s)ds$$
where, $\varphi^N(s)=\Lambda(0)X^N_s(1-X^N_s).$
\end{lemma}
\bpf
For each $1\le i<N$, let $P^{i}$ be a Poisson process with intensity $\Lambda(0)(N-i)$. At time $t\in P^i$, we have 
 
 \begin{eqnarray*}
\Delta X^N_t=
\begin{cases}
\frac{1}{N}, &  \mbox{if $\eta_{t^-}(i)=1$ and $\eta_{t^-}(N)=0$}  \\
-\frac{1}{N}, & \mbox{if $\eta_{t^-}(i)=0$ and $\eta_{t^-}(N)=1$} \\
0 , &     \mbox{otherwise.}
\end{cases}
\end{eqnarray*}
 Now, let 
\begin{align*}
 A_i&= \{\eta_{t}(i)=1,\eta_{t}(N)=0\},\\
 B_i&=\{\eta_{t}(i)=0,\eta_{t}(N)=1\}.\\
 \end{align*}
 We have 
 $$\P(A_i\mid X^N_t)=\P(B_i\mid X^N_t)=\frac{N}{N-1}X^N_t(1-X^N_t),$$
from which, we deduce that
\begin{align*}
\langle \mathcal{K}^N \rangle_t&=\frac{1}{N^2}\Lambda(0)\frac{N(N-1)}{2}\frac{2N}{N-1}X^N_t(1-X^N_t)\\
&=\Lambda(0)X^N_t(1-X^N_t)
\end{align*}
The result is proved. 

\epf

\bigskip

Now, let 
\begin{align*}
Y^N_t&=X^N_0+\int_{ [0, t] \times]0, 1]^4 }\psi^N(X^N_{s^-}, u, p, v, w)\bar{M}(ds, du, dp, dv, dw)\\
&-\frac{1}{N}\int_{[0, t]\times[0, 1]}{\bf 1}_{u\le X^N_{s^-}}{\bf 1}_{\eta_{s-}(N+1)=0}M^N_1(ds, du).
\end{align*}
We have 
\begin{equation}\label{NXN2}
X^N_t=\mathcal{K}^N_t+Y^N_t, \quad \forall t\ge 0
\end{equation}
From lemma \ref{crochM}, we have $\forall T\ge 0$
$$\sup_{0\le t\le T}\sup_{N\ge 1}\mid\varphi^N(s)\mid\leq C\quad a.s.$$
Using the last identity, we deduce by Aldous' tightness criterion (see Aldous \cite{A.T}) that 
$$\{\mathcal{K}^N_t, t\ge 0, N\ge 1\} \; \mbox{is tight in}\  D([0, \infty)) .$$
Since $\mathcal{K}^N$ is tight, there exists a subsequence of the sequence $\mathcal{K}^N$ such that 
$$\mathcal{K}^N\Rightarrow \mathcal{K}\ \mbox{weakly in} \; D([0, \infty)),$$
where $\mathcal{K}$ is a continuous martingale (since the jumps of $\mathcal{K}^N$ are of size $\pm\frac{1}{N}$) such that
\begin{equation}\label{crochl}
<\mathcal{K}>_t=\int_{0}^{t}cX_s(1-X_s)ds,
\end{equation}
where $c=\Lambda (0)$. The main result of this section is 

\begin{theorem}\label{main-result2}
Suppose that $X^N_0\rightarrow x$ a.s, as $N\rightarrow \infty$. Then the $[0,1]-$valued process $\{X_t ,t \ge 0\}$ defined by \eqref{NX-limit} is the (unique in law) solution to the  stochastic differential equation
\begin{align}\label{KLWF}
\begin{split}
X_t = x &- \alpha \int_{0}^{t} X_s(1-X_s)ds\\
&+\int_{0}^{t}\sqrt {\Lambda (0)X_s(1-X_s)}dB_s\\
&+\int_{[0, t]\times]0, 1[^2}p({\bf 1}_{u\le X_{s^-}}-X_{s^-})\bar{M}(ds, du, dp),
\end{split}
\end{align}
where $\bar{M}$ is the compensated measure $M$ defined in section 3.2, and  
$B$ is a standard Brownian motion.
\end{theorem}
The identification of the limiting equation is done similarly as in the proof of Theorem \ref{main-result}. Strong uniqueness of the solution to \eqref{KLWF} follows again from Dawson and Li \cite{DL}, and weak uniqueness could also be proved by a duality argument.

Since Kingman's coalescent comes down from infinity, we have fixation in our new model in finite time as soon as $\Lambda(0)>0$. 

\subsection*{Acknowledgements}
We thank an anonymous Referee, who drew our attention to the work of H\'enard \cite{Henard}, which
permitted us to add the subsection \ref{ss:hn} to our original version. Other Referee's remarks helped us 
to improve other points of our first version.

{\small
This work has been supported by the ANR project MANEGE, and by the Infectiop\^ole Sud Foundation.
}

Boubacar Bah I2M, Aix--Marseille Universit\'e, 39, rue F. Joliot Curie, F 13453 Marseille cedex 13. bbah12@yahoo.fr

Etienne Pardoux (corresponding author) Aix-Marseille Universit\'e, CNRS, Centrale Marseille, I2M, UMR 7373 13453
Marseille, France. 
etienne.pardoux@univ-amu.fr

\end{document}